\documentclass[11pt]{article}
\begin{document}
\begin{center}{\bf On boundedness, existence and uniqueness of strong solutions of the Navier-Stokes Equations in 3 dimenions.}

A. A. Ruzmaikina.

{\it Department of Statistics, Purdue University, West Lafayette, IN 47909.}

{\tt aar@purdue.edu}
\end{center}

\begin{center}{\bf Abstract.}\end{center}
{\it In this paper we consider the Navier-Stokes Equations in ${\bf R}^3$ in the vorticity formulation in the absence of the external forces. We first derive the upper bounds on $|\omega|_{\frac{6.02}{3}}(t)$ for $t \in [t_1, t_2]$ in terms of ${\rm max}_{t \in [t_1, t_2]} |\omega|_2(t)$, then derive the upper bounds on $|{\bar \nabla}{\bar \omega}|_2(t)$ for $t \in [t_1, t_2]$ in terms of ${\rm max}_{t \in [t_1, t_2]} |\omega|_2(t)$ and $|\omega|_{\frac{6.02}{3}}(t_1)$; then derive the upper bound on $|{\bar \nabla}{\bar \omega}|_{3.01}(t)$ for $t \in [t_1, t_2]$ in terms of ${\rm max}_{t \in [t_1, t_2]} |\omega|_2(t)$, $|\omega|_{\frac{6.02}{3}}(t_1)$ and $|{\bar \nabla}{\bar \omega}|_2(t_1)$; then we derive the upper bound on $|\omega|_{\infty}(t)$ for $t \in [t_1, t_2]$ in terms of ${\rm max}_{t \in [t_1, t_2]} |\omega|_2(t)$, $|\omega|_{\frac{6.02}{3}}(t_1)$, $|{\bar \nabla}{\bar \omega}|_2(t_1)$ and ${\bar \nabla}{\bar \omega}|_{3.01}(t_1)$.
These estimates are used to improve the estimate on $|\omega|_{\infty}(t)$ and to show that $|\omega|_{\infty}(t) \leq {\rm max}(|\omega|_{\infty}(t_1), e^e)$, provided that certain assumptions hold. Repeated use of Local Existence and Uniqueness results together with the upper bounds on $|\omega|_{\infty}(t)$ gives the existence and uniqueness of the solution if $|\omega|_{\infty}(0) < \infty$ or if $|\omega|_4(0) < \infty.$}

\begin{center}{\bf Introduction and statement of results.} \end{center}

The Navier-Stokes Equations for velocity can be written as
\begin{eqnarray}
&& \partial_t {\bar u}({\bar x}, t) - \nu \Delta {\bar u}({\bar x}, t) +({\bar u}({\bar x}, t) \cdot {\bar \nabla}){\bar u}({\bar x}, t) + {\bar \nabla} p({\bar x}, t) = f({\bar x}, t), \\ \nonumber && \qquad \qquad \qquad {\bar \nabla} \cdot {\bar u}({\bar x}, t) = 0, ~ {\bar x} \in {\bf R}^3,~ t \in [0, \infty).
\end{eqnarray}
We assume that $f({\bar x}, t) =0$ and that ${\rm lim}_{|{\bar x}| \rightarrow \infty} |{\bar u}({\bar x}, t)| = 0.$ 

Introducing vorticity ${\bar \omega}({\bar x}, t) = {\bar \nabla} \times {\bar u}({\bar x}, t)$ allows us to obtain the Navier-Stokes Equations without the pressure term.

The Navier-Stokes Equations for vorticity can be written as
\begin{eqnarray}
&& \partial_t {\bar \omega}({\bar x}, t) - \nu \Delta {\bar \omega}({\bar x}, t) = -({\bar u}({\bar x}, t) \cdot {\bar \nabla}) {\bar \omega}({\bar x}, t) + ({\bar \omega}({\bar x}, t) \cdot {\bar \nabla}) {\bar u}({\bar x}, t).\ \\
\nonumber && \qquad \qquad \qquad {\bar \nabla} \cdot {\bar \omega} ({\bar x}, t) = 0, ~ {\bar x} \in {\bf R}^3,~ t \in [0, \infty).
\end{eqnarray}
We assume that ${\rm lim}_{|x| \rightarrow \infty} |{\bar \omega}({\bar x}, t)| = 0.$

The question of existence and uniqueness of a strong solution of these equations in 3 dimensions is a longstanding problem. The solution of this problem is described in this paper.

The main result of the paper are the following Theorems:

{\it {\bf Theorem 1:} Suppose that $|\omega|_2 (t) < \infty$ for $t \in [t_1, t_2]$ and also that $|\omega|_{\infty}(t_1) < \infty$ and  $|{\bar \nabla}{\bar \omega}|_{3.01} (t_1) < \infty$. Then $|\omega|_{\infty}(t) < \infty$ for $t \in [t_1, t_2]$ and an upper bound (equation (\ref{omegainftyupper})) on $|\omega|_{\infty}(t)$ can be given as a function of $t$, ${\rm max}_{t \in [t_1, t_2]} |\omega|_2(t)$, $|\omega|_{\infty}(t_1)$ and $|{\bar \nabla}{\bar \omega}|_{3.01}(t_1)$.}

{\it {\bf Theorem 2:} Suppose that $|\omega|_2 (t) < \infty$ for $t \in [t_1, t_2]$ and that $|\omega|_{\infty}(t_1) < \infty$ and $|{\bar \nabla}{\bar \omega}|_{3.01} (t_1) < \infty$. Suppose also that there exists $n_0$ such that for all  $n > n_0$, either $|\omega|_n(t) \geq e^e$ for all $t \in [t_1, t_2]$ or $|\omega|_n(t)$ is a continuous function of $t$ for all $t \in [t_1, t_2].$
Then for all sufficiently large $n$,  $|\omega|_n(t) \leq {\rm max}(e^{\frac{1}{n}e^{e^{t_2-t_1}}}|\omega|_n(t_1), e^e)$ and $|\omega|_{\infty}(t) \leq {\rm max}(|\omega|_{\infty}(t_1), e^e)$ for $t \in [t_1, t_2]$.}

{\it {\bf Theorem 3:} Suppose that $|\omega|_{\infty}(0) < \infty$. Then for any $T > 0$, there exists a unique solution  of the 3D Navier Stokes Equations ${\bar \omega}({\bar x},t)$ for $x \in {\bf R}^3$ and $t \in [0, T]$, such that $|\omega|_{\infty}(t) < \infty$ for $t \in [0, T]$.}

{\it {\bf Theorem 4:} Suppose that $|\omega|_4 (0) < \infty$. Then for any $T > 0$, there exists a unique solution of the 3D Navier Stokes Equations ${\bar \omega}({\bar x}, t)$ for $x \in {\bf R}^3$ and $t \in [0, T]$, such that $|\omega|_{\infty}(t) < \infty$ for $t \in (0, T]$.}

\begin{center} {\bf Proof of Theorem 1.}\end{center}

{\bf Part 1:} Take Navier Stokes equation for ${\bar \omega}$ and take the inner product with ${\bar \omega}(x)$ integrate over space and by parts on the viscosity term to obtain  
\begin{eqnarray}\label{nsomega2}
\frac{1}{2}\partial_t \int |{\bar \omega}({\bar x})|^2 d^3 x + \nu \int |{\bar \nabla}{\bar \omega}({\bar x})|^2 d^3 x = \frac{3}{4\pi}\int d^3 x \int d^3 y ({\bar \omega}({\bar x}) \cdot {\hat y}) ({\bar \omega}({\bar x}) \times {\bar \omega}({\bar x}+{\bar y}) \cdot {\hat y}) \frac{1}{|y|^3}.
\end{eqnarray}
The absolute value of the integral on the right hand side of (\ref{nsomega2}) can be estimated as
\begin{eqnarray}
&& \left|\int d^3x \int d^3y ({\bar \omega}({\bar x}) \cdot {\hat y}) ({\bar \omega}({\bar x}) \times {\bar \omega}({\bar x}+{\bar y}) \cdot {\hat y}) \frac{1}{|y|^3}\right| \leq \\ \nonumber && \leq
\left| \int d^3 x \int_{|y| < 1} ({\bar \omega}({\bar x}) \cdot {\hat y}) ({\bar \omega}({\bar x}) \times ({\bar \omega}({\bar x}+{\bar y})-{\bar \omega}({\bar x})) \cdot {\hat y}) \frac{d^3y}{|y|^3}\right| + \\ \nonumber && \qquad \qquad \qquad + \left| \int d^3 x \int_{|y| \geq 1} ({\bar \omega}({\bar x}) \cdot {\hat y}) ({\bar \omega}({\bar x}) \times {\bar \omega}({\bar x}+{\bar y}) \cdot {\hat y}) \frac{d^3 y}{|y|^3} \right| \\ \nonumber && \leq
3 \sqrt{6}\int_0^1 ds \int d^3x \int_{|y| < 1} |{\bar \omega}({\bar x})|^2 |{\bar \nabla}{\bar \omega}({\bar x}+s{\bar y})| \frac{d^3y}{|y|^2}  \\ \nonumber && \quad \quad \quad +3 \sqrt{2} \int d^3x |{\bar \omega}({\bar x})|^2 \int_{|y| \geq 1} \frac{|{\bar \omega}({\bar x}+{\bar y})|}{|y|^3} d^3y \\ \nonumber && \leq 3 \sqrt{6}\int_0^1 ds \int d^3 x' |{\bar \nabla}{\bar \omega}({\bar x'})| \int_{|y| < 1} d^3y  \frac{|{\bar \omega}({\bar x'}-s{\bar y})|^2}{|y|^2} +\\ \nonumber && \quad \quad \quad + 3\sqrt{2}\int d^3x |{\bar \omega}({\bar x})|^2 \int_{|y| \geq 1} d^3y \frac{|{\bar \omega}({\bar x} +{\bar y})|}{|y|^3}
\\ \nonumber &&\leq 3 \sqrt{6} \int_0^1 ds \sqrt{\int d^3x' |{\bar \nabla}{\bar \omega}({\bar x'})|^2} \times \\ \nonumber && 
 \sqrt{\int d^3x' \int_{|y| < 1 } \frac{|{\bar \omega}({\bar x'} - s{\bar y})|^2}{|y|^2} d^3y \int_{|z| < 1} \frac{|{\bar \omega}({\bar x'}-s{\bar z})|^2}{|z|^2}d^3z }  + \sqrt{24 \pi}\int d^3x |{\bar \omega}({\bar x})|^2 (\int d^3y |{\bar \omega}({\bar y})|^2)^{\frac{1}{2}} \\ \nonumber && 
\leq 3\sqrt{6} \int_0^1 ds \sqrt{\int d^3x |{\bar \nabla}{\bar \omega}({\bar x})|^2} \sqrt{\int d^3x' \int_{|y| \leq 1} \frac{|{\bar \omega}({\bar x'}-s{\bar y})|^2}{|y|^{3-\epsilon}} d^3y  \int_{|z| \leq 1} \frac{|{\bar \omega}({\bar x'}-s{\bar z})|^2}{|z|^{1+\epsilon}} d^3z}+\\ \nonumber && \quad \quad \quad + \sqrt{24 \pi} |\omega|_2^2 |\omega|_2 \\ \nonumber &&  \leq 
3 \sqrt{6} \int_0^1 ds \sqrt{\int d^3x |{\bar \nabla}{\bar \omega}({\bar x})|^2} \sqrt{\int d^3 x |{\bar \omega}({\bar x})|^{2} \int_{|y| \leq 1} \frac{d^3 y}{|y|^{3-\epsilon}} \int_{|z| \leq 1} \frac{|{\bar \omega}({\bar x}+s{\bar y}-s{\bar z})|^2}{|z|^{1+\epsilon}} d^3z}  +\\ \nonumber && \quad \quad \quad +  \sqrt{24 \pi}|\omega |_2^2 |\omega |_2\\ \nonumber && \leq 3 \sqrt{6} \int_0^1 \frac{ ds}{s^{1-\frac{\epsilon}{2}}} \sqrt{\int d^3x |{\bar \nabla}{\bar \omega}({\bar x})|^2 } \sqrt{\int d^3 x |{\bar \omega}({\bar x})|^2 \int_{|y| \leq 1} \frac{d^3y}{|y|^{3-\epsilon}} \int_{|z'| \leq 1} \frac{|{\bar \omega}({\bar x}+s{\bar y}-{\bar z'})|^2}{|z'|^{1+\epsilon}} d^3z'} +\\ \nonumber && \quad \quad \quad +  \sqrt{24 \pi} |\omega|_2^2 |\omega|_2 \\ \nonumber && \leq 3 \sqrt{6} \int_0^1 \frac{ds}{s^{1-\frac{\epsilon}{2}}} \sqrt{\int d^3x |{\bar \nabla}{\bar \omega}({\bar x})|^2}\\ \nonumber &&  \sqrt{\int d^3x |{\bar \omega}({\bar x})|^2 \int_{|y| \leq 1} \frac{d^3 y}{|y|^{3-\epsilon}} (\int d^3 z'|{\bar \omega}({\bar x} - {\bar z'} + s{\bar y})|^{\frac{2(3-3\epsilon)}{2-3\epsilon}})^{\frac{2-3\epsilon}{3-3\epsilon}} (\int_{|z'| < 1}\frac{d^3 z'}{|z'|^{3(1-\epsilon^2)}})^{\frac{1}{3(1-\epsilon)}} }+\\ \nonumber && \quad \quad \quad +\sqrt{24 \pi}|\omega |_2^2 |\omega |_2 \\ \nonumber && \leq 
3 \sqrt{6}(\frac{4 \pi}{\epsilon})^{\frac{1}{2}}(\frac{4\pi}{3\epsilon^2})^{\frac{1}{6(1-\epsilon)}}\int_0^1 \frac{ds}{s^{1-\frac{\epsilon}{2}}} \sqrt{\int d^3x |{\bar \nabla}{\bar \omega}({\bar x})|^2} \sqrt{\int d^3x |{\bar \omega}({\bar x})|^2(\int d^3z |{\bar \omega}({\bar z})|^{\frac{2(3-3\epsilon)}{2-3\epsilon}})^{\frac{2-3\epsilon}{3-3\epsilon}}}\\ \nonumber && \quad \quad \quad  + \sqrt{24 \pi}|\omega |_2^2 |\omega |_2 \\ \nonumber && \leq \frac{ 6\sqrt{6}}{\epsilon}(\frac{4\pi}{\epsilon})^{\frac{1}{2}}(\frac{4 \pi}{3\epsilon^2})^{\frac{1}{6(1-\epsilon)}} (\int |{\bar \nabla}{\bar \omega}({\bar x})|^2 d^3x )^{\frac{3-2\epsilon}{4(1-\epsilon)}} (\int |{\bar \omega}({\bar x})|^2 d^3x)^{\frac{3-4\epsilon}{4-4\epsilon}} +  \sqrt{24 \pi}|\omega |_2^2 |\omega |_2 \\ \nonumber && \leq \frac{1-2\epsilon}{4-4\epsilon}(\frac{6 \sqrt{6}}{\epsilon})^{\frac{4-4\epsilon}{1-2\epsilon}} (\frac{4 \pi}{\epsilon})^{\frac{2-2\epsilon}{1-2\epsilon}}(\frac{4 \pi}{3 \epsilon^2})^{\frac{2}{3(1-2\epsilon)}} (\frac{2^{0.33}}{\nu})^{\frac{3-2\epsilon}{1-2\epsilon}} (\int d^3x |{\bar \omega}({\bar x})|^2)^{\frac{3-4\epsilon}{1-2\epsilon}}+\\ \nonumber &&\quad \quad \quad + \frac{\nu}{2^{0.33}} \int d^3x |{\bar \nabla}{\bar \omega}({\bar x})|^2+ \sqrt{24 \pi}|\omega |_2^2 |\omega |_2.
\end{eqnarray}
Notice that in the above estimate $\epsilon$ is arbitrarily small.
Using this estimate and introducing $\varepsilon$ such that $$\varepsilon = \frac{2\epsilon}{1-2\epsilon}$$ we obtain the following inequality from (\ref{nsomega2})
\begin{eqnarray}
&& \frac{1}{2} \partial_t \int |{\bar \omega}({\bar x})|^2 d^3x + \nu(1-\frac{1}{2^{0.33}}) \int |{\bar \nabla}{\bar \omega}({\bar x})|^2 d^3x \leq \nonumber \\ && \qquad \qquad \qquad \leq \frac{3}{4 \pi}\left(\frac{0.331 \cdot 6^{6}\cdot (4 \pi)^{\frac{8}{3}}}{\varepsilon^{\frac{22}{3}+\frac{13}{3} \varepsilon} \nu^{3+2 \varepsilon}} (\int |{\bar \omega}({\bar x})|^2 d^3x)^{3+\varepsilon}+  \sqrt{24 \pi}|\omega|_2^3\right).\qquad \qquad \qquad 
\end{eqnarray}
Note that $\varepsilon$ can be chosen arbitrarily small. Below we assume that $\varepsilon = \frac{1}{100}.$

Since $\int |{\bar \omega}({\bar x})|^2 d^3 x < \infty$ for $t \in [t_1, t_2]$, there exists $N > 0$ such that  $\int |{\bar \omega}({\bar x})|^2 d^3 x (t) \leq N$ for $t \in [t_1, t_2]$. Denote $$C(t_1, t_2) = \int_{t_1}^{t_2} dt \int |{\bar \omega}({\bar x})|^2 d^3 x, C(t) = \int_0^t ds \int |{\bar \omega}({\bar x})|^2 d^3x (s).$$
Then we obtain 
\begin{eqnarray}
&&\frac{1}{2} \int |{\bar \omega}({\bar x})|^2 d^3 x (t_2) - \frac{1}{2}\int |{\bar \omega}({\bar x})|^2 d^3 x (t_1) + \nu(1-\frac{1}{2^{0.33}}) \int_{t_1}^{t_2} dt \int |{\bar \nabla}{\bar \omega}({\bar x})|^2 d^3 x\nonumber \\ \nonumber &&  \leq \frac{2.74 \cdot 10^{6}}{\varepsilon^{\frac{22}{3}+\frac{13}{3}\varepsilon} \nu^{3+2\varepsilon}} N^{2+\varepsilon}C(t_1,t_2)+2.1 N^{\frac{1}{2}} C(t_1,t_2) \leq \frac{2.35 \cdot 10^{21}}{\nu^{3.02}} N^{2+\varepsilon}C(t_1, t_2)+2.1 N^{\frac{1}{2}} C(t_1, t_2).
\end{eqnarray}
Therefore
\begin{eqnarray}
 && \int_{t_1}^{t_2} dt \int d^3 x |{\bar \nabla}{\bar \omega}({\bar x})|^2 \leq \frac{1.075 \cdot 10^{22}}{ \nu^{4+2\varepsilon}} N^{2+\varepsilon} C(t_1, t_2) +\frac{9.6}{\nu} N^{\frac{1}{2}} C(t_1, t_2)+  \frac{4.57}{\nu} N. \qquad
\end{eqnarray}

{\bf Part 2:} Take the inner product of the Navier-Stokes equation for ${\bar \omega}({\bar x})$ with $|\omega({\bar x})|^{\frac{2}{3}\varepsilon} {\bar \omega}({\bar x})$, integrate over space and integrate the viscosity term by parts to obtain
\begin{eqnarray}\label{nsom23}
&& \frac{1}{2+\frac{2}{3}\varepsilon} \partial_t \int |\omega|^{2+\frac{2}{3}\varepsilon} ({\bar x}) d^3 x + \nu \int {\bar \nabla}(|\omega|^{\frac{2}{3}\varepsilon} ({\bar x}) {\bar \omega}({\bar x})) ({\bar \nabla} {\bar \omega}({\bar x})) d^3 x  \\ \nonumber && =
\frac{3}{4 \pi} \int \int (|\omega|^{\frac{2}{3}\varepsilon} ({\bar x}) {\bar \omega}({\bar x}) \cdot {\hat y}) ({\bar \omega}({\bar x}) \times {\bar \omega}({\bar x} + {\bar y}) \cdot {\hat y}) \frac{d^3 x d^3 y}{|y|^3}.
\end{eqnarray}
The absolute value of the integral on the right hand side of (\ref{nsom23}) can be estimated from above as
\begin{eqnarray}
&& \left|\int \int (|\omega|^{\frac{2}{3}\varepsilon}({\bar x}) {\bar \omega}({\bar x}) \cdot {\hat y}) ({\bar \omega}({\bar x}) \times {\bar \omega}({\bar x} + {\bar y}) \cdot {\hat y}) \frac{d^3 x d^3 y}{|y|^3}\right| \leq \\ \nonumber && \leq 
\left|\int \int_{|y| < 1} (|\omega|^{\frac{2}{3}\varepsilon}({\bar x}) {\bar \omega}({\bar x}) \cdot {\hat y}) ({\bar \omega}({\bar x}) \times ({\bar \omega}({\bar x} + {\bar y}) - {\bar \omega}({\bar x})) \cdot {\hat y}) \frac{d^3 x d^3 y}{|y|^3}\right| + \\ \nonumber && + \left| \int \int_{|y| \geq 1} (|\omega|^{\frac{2}{3}\varepsilon}({\bar x}){\bar \omega}({\bar x}) \cdot {\hat y}) ({\bar \omega}({\bar x}) \times {\bar \omega}({\bar x} + {\bar y}) \cdot {\hat y}) \frac{d^3 x d^3 y}{|y|^3}\right| \leq \\ \nonumber && \leq
3 \sqrt{6} \int_{0}^{1} ds \int d^3x  |\omega|^{2+\frac{2}{3}\varepsilon} ({\bar x}) \int_{|y| <1 } \frac{|{\bar \nabla} {\bar \omega}({\bar x} + s {\bar y})|}{|y|^2} d^3 y+ \\ \nonumber && \qquad \qquad \qquad \qquad \qquad  + 3 \sqrt{2}\int \int_{|y| \geq 1} |\omega|^{2+\frac{2}{3}\varepsilon}({\bar x}) |\omega|({\bar x} + {\bar y}) \frac{d^3 x d^3 y}{|y|^3} \leq \\ \nonumber && \leq 3 \sqrt{6}\int_{0}^{1} ds \int d^3 x' |{\bar \nabla}{\bar \omega}({\bar x'})| \int_{|y| < 1} \frac{|\omega|^{2+\frac{2}{3}\varepsilon}({\bar x'} - s {\bar y})}{|y|^2} d^3 y + \\ \nonumber && \qquad \qquad \qquad \qquad \qquad + 3 \sqrt{2} \int |\omega|^{2+\frac{2}{3}\varepsilon}({\bar x}) d^3 x \int_{|y| \geq 1} \frac{|{\bar \omega}({\bar x}+{\bar y})|}{|y|^3} d^3 y \leq \\ \nonumber && \leq 3 \sqrt{6} \int_{0}^{1} ds  \sqrt{\int d^3 x' \int_{|y| < 1} \frac{|\omega|^{2+\frac{2}{3}\varepsilon}({\bar x'} - s {\bar y})}{|y|^2} d^3 y \int_{|z| < 1} \frac{|\omega|^{2+\frac{2}{3}\varepsilon}({\bar x'} - s {\bar z})}{|z|^2} d^3 z} |{\bar \nabla}{\bar \omega}|_2   \\ \nonumber && \qquad \qquad \qquad \qquad \qquad 
+ \sqrt{24 \pi} |{\omega}|_{2+\frac{2}{3}\varepsilon}^{2+\frac{2}{3}\varepsilon} |{\omega}|_2 \\ \nonumber && \leq 3 \sqrt{6}\int_{0}^{1} ds \sqrt{\int d^3 x' \int_{|y| \leq 1} \frac{ |\omega|^{2 + \frac{2}{3} \varepsilon}({\bar x'}-s{\bar y})}{|y|^{3-\varepsilon}} d^3 y \int_{|z| \leq 1} \frac{|\omega|^{2+\frac{2}{3}\varepsilon}({\bar x'} - s{\bar z})}{|z|^{1+\varepsilon}} d^3 z} \sqrt{\int |{\bar \nabla} {\bar \omega} ({\bar x})|^2 d^3 x}\\ \nonumber && \qquad \qquad \qquad \qquad \qquad  +\sqrt{24 \pi} |\omega|_{2+\frac{2}{3}\varepsilon}^{2+\frac{2}{3}\varepsilon} |\omega|_2 \\ \nonumber && \leq 3 \sqrt{6} \int_{0}^{1} ds \sqrt{\int |\omega|^{2 + \frac{2}{3} \varepsilon}({\bar x''}) d^3 x''\int_{|y| \leq 1} \frac{d^3 y}{|y|^{3-\varepsilon}} \int_{|z| \leq 1} \frac{|\omega|^{2+\frac{2}{3}\varepsilon}({\bar x''}-s({\bar z} - {\bar y}))}{|z|^{1+\varepsilon}} d^3 z} \ |{\bar \nabla}{\bar \omega}|_2 \\ \nonumber && \qquad \qquad \qquad \qquad \qquad  + \sqrt{24 \pi}|\omega|_{2+\frac{2}{3}\varepsilon}^{2+\frac{2}{3}\varepsilon} |\omega|_2 \leq \\ \nonumber && \leq 3 \sqrt{6} \int_{0}^{1} ds \sqrt{\int d^3 x'' |\omega|^{2+\frac{2}{3}\varepsilon} (x'') \int_{|y| \leq 1} \frac{d^3 y}{|y|^{3-\varepsilon}} \int_{|z'| \leq s} \frac{|\omega|^{2+\frac{2}{3}\varepsilon}(x''+sy - z')}{|z'|^{1+\varepsilon}} \frac{d^3 z'}{s^{2-\varepsilon}}}|{\bar \nabla} {\bar \omega}|_2\\ \nonumber && \qquad \qquad \qquad \qquad \qquad + \sqrt{24\pi}|\omega|_{2+\frac{2}{3}\varepsilon}^{2+\frac{2}{3}\varepsilon}|\omega|_2 \leq \\ \nonumber && \leq 3 \sqrt{6}\int_{0}^{1} \frac{ds}{s^{1-\frac{\varepsilon}{2}}} \sqrt{\int d^3 x |\omega|^{2+\frac{2}{3}\varepsilon}({\bar x}) \int_{|y| \leq 1} \frac{d^3 y}{|y|^{3-\varepsilon}} \int_{|z'| \leq 1} \frac{|\omega|^{2+\frac{2}{3} \varepsilon}(x+sy-z')}{|z'|^{1+\varepsilon}}d^3z'}|{\bar \nabla}{\bar \omega}|_2+ \\ \nonumber &&\qquad \qquad \qquad \qquad \qquad   +\sqrt{24 \pi} |\omega|_{2+\frac{2}{3}\varepsilon}^{2+\frac{2}{3}\varepsilon} |\omega|_2 \leq \\ \nonumber && \leq 3 \sqrt{6} (\frac{4 \pi}{3\varepsilon^2})^{\frac{1}{6(1-\varepsilon)}}\\ \nonumber && \ \ \ \int_0^1 \frac{ds}{s^{1-\frac{\varepsilon}{2}}} \sqrt{\int d^3 x |\omega|^{2+\frac{2}{3}\varepsilon} ({\bar x}) \int_{|y| \leq 1} \frac{d^3 y}{|y|^{3-\varepsilon}} (\int|\omega|^{(2+\frac{2}{3}\varepsilon) \frac{3-3\varepsilon}{2-3\varepsilon}} (x+sy - z')d^3 z')^{\frac{2-3\varepsilon}{3-3\varepsilon}}}|{\bar \nabla}{\bar \omega}|_2  \\ \nonumber && \qquad \qquad \qquad \qquad \qquad + \sqrt{24 \pi} |\omega|_{2+\frac{2}{3}\varepsilon}^{2+\frac{2}{3}\varepsilon} |\omega|_2 \\ \nonumber && \leq 
3 \sqrt{6} (\frac{4 \pi}{\varepsilon})^{\frac{1}{2}}(\frac{4\pi}{3\varepsilon^2})^{\frac{1}{6(1-\varepsilon)}}|{\bar \nabla}{\bar \omega}|_2 \\ \nonumber && \int_{0}^{1} \frac{ds}{s^{1-\frac{\varepsilon}{2}}} \sqrt{(\int |{\bar \omega}({\bar x})|^{2+\frac{2}{3} \varepsilon}d^3 x)^{1+\frac{(1-\frac{7}{3}\varepsilon+\frac{\varepsilon^2}{3})(1-\frac{3}{2}\varepsilon)}{2(1-\frac{5}{3}\varepsilon+\frac{\varepsilon^2}{4})(1-\varepsilon)}} (\int |{\bar \nabla}{\bar \omega}({\bar x})|^2 d^3 x)^{\frac{(1+\frac{\varepsilon}{3})(1-\frac{3}{2}\varepsilon)}{2(1-\frac{5}{3}\varepsilon+\frac{\varepsilon^2}{4})(1-\varepsilon)}}} + \\ \nonumber && \qquad \qquad \qquad \qquad \qquad + \sqrt{24 \pi} |\omega|_{2+\frac{2}{3}\varepsilon}^{2+\frac{2}{3}\varepsilon} |\omega|_2 \\ \nonumber && \leq \frac{6\sqrt{6}}{\varepsilon}(\frac{4 \pi}{\varepsilon})^{\frac{1}{2}}(\frac{4\pi}{3\varepsilon^2})^{\frac{1}{6(1-\varepsilon)}}(\int|\omega|^{2+\frac{2}{3}\varepsilon}({\bar x}) d^3x)^{\frac{1}{2}+\frac{(1 - \frac{7}{3} \varepsilon+\frac{\varepsilon^2}{3})(1-\frac{3}{2}\varepsilon)}{4(1-\frac{5}{3}\varepsilon+\frac{\varepsilon^2}{4})(1-\varepsilon)}}  (\int|{\bar \nabla} {\bar \omega}({\bar x})|^2 d^3 x)^{\frac{1}{2}+\frac{(1+\frac{\varepsilon}{3})(1-
\frac{3}{2}\varepsilon)}{4(1-\frac{5}{3}\varepsilon+\frac{\varepsilon^2}{4})(1-\varepsilon)}} \\ \nonumber && \qquad \qquad \qquad \qquad \qquad +\sqrt{24 \pi} |\omega|_{2+\frac{2}{3}\varepsilon}^{2+\frac{2}{3}\varepsilon} |\omega|_2 \\ \nonumber && \leq 3.0618 \cdot 10^5 \left(\int |{\bar \nabla}{\bar \omega}({\bar x})|^2 d^3 x+1\right) \left(\int |{\bar \omega}({\bar x})|^{2 + \frac{2}{3}\varepsilon} d^3 x\right)^{\alpha}+  \sqrt{24 \pi} |\omega|_{2+\frac{2}{3}\varepsilon}^{2+\frac{2}{3}\varepsilon} |\omega|_2,
\end{eqnarray} 
where $\alpha = \frac{1}{2} + \frac{1}{4}\frac{(1-\frac{7}{3}\varepsilon+\frac{\varepsilon^2}{3})(1-\frac{3}{2}\varepsilon)}{(1-\frac{5}{3}\varepsilon+\frac{\varepsilon^2}{4})(1-\varepsilon)}\approx 0.747.$
Putting the estimates together we obtain
\begin{eqnarray}\label{rhsns223}
&& \frac{1}{2+\frac{2}{3}\varepsilon} \partial_t \int |{\bar \omega}({\bar x})|^{2+\frac{2}{3} \varepsilon} d^3 x + \nu \int {\bar \nabla}(|{\bar \omega}({\bar x})|^{\frac{2}{3}\varepsilon}{\bar \omega}({\bar x}))({\bar \nabla}{\bar \omega}({\bar x})) d^3 x \leq \qquad \qquad \qquad \\ \nonumber &&\quad \quad \quad \leq  \frac{3}{4\pi} 3.0618 \cdot 10^5 \left(\int |{\bar \omega}({\bar x})|^{2+\frac{2}{3}\varepsilon}d^3 x\right)^{\alpha} \left(\int |{\bar \nabla}{\bar \omega}({\bar x})|^2 d^3 x+1 \right)+\frac{3 \sqrt{24 \pi}}{4 \pi} |\omega|_{2+\frac{2}{3}\varepsilon}^{2+\frac{2}{3}\varepsilon} N^{\frac{1}{2}}\\ \nonumber && \quad \quad \quad \leq 0.7305 \cdot 10^5 \left(\int |\omega|^{2+\frac{2}{3}\varepsilon} d^3 x\right)^{\alpha} \left( \int |{\bar \nabla}{\bar \omega}({\bar x})|^2 d^3 x + 1\right) + 2.1 |\omega|_{2+\frac{2}{3}\varepsilon}^{2+\frac{2}{3}\varepsilon} N^{\frac{1}{2}}.
\end{eqnarray}
Since we can rewrite the viscosity term as shown in the ({\ref{viscosity223}) we obtain that the viscosity term is positive definite
\begin{eqnarray}\label{viscosity223}
\int {\bar \nabla} (|\omega|^{\frac{2}{3}\varepsilon}({\bar x}) {\bar \omega}({\bar x})) ({\bar \nabla} {\bar \omega}) d^3 x = \frac{\frac{2}{3}\varepsilon}{(1+\frac{\varepsilon}{3})^2} 
\int ({\bar \nabla}|\omega|^{1+\frac{\varepsilon}{3}} ({\bar x}))^2 d^3 x + \int |\omega|^{\frac{2}{3}\varepsilon} ({\bar x}) ({\bar \nabla} {\bar \omega} ({\bar x}))^2 d^3 x \qquad
\end{eqnarray}
From (\ref{rhsns223}) we obtain
\begin{eqnarray}
&& \int |{\bar \omega}({\bar x})|^{2+\frac{2}{3}\varepsilon} d^3 x (t_2) \leq \\ \nonumber && \leq e^{4.22 \ N^{\frac{1}{2}}(t_2-t_1)} \left(0.371 \cdot 10^5 \int_{t_1}^{t_2} dt \left(\int |{\bar \nabla}{\bar \omega}|^2({\bar x}) d^3 x+1\right)+ \left(\int |{\bar \omega}({\bar x})|^{2+\frac{2}{3}\varepsilon} d^3 x (t_1)\right)^{1-\alpha}\right)^{\frac{1}{1-\alpha}} \\ \nonumber && \leq e^{4.22 N^{\frac{1}{2}}(t_2 - t_1)} \\ \nonumber && \left(\frac{3.99 \cdot 10^{26}}{\nu^{4+2\varepsilon}} N^{2+\varepsilon} C(t_1, t_2) + \frac{3.562 \cdot 10^5}{\nu} N^{\frac{1}{2}} C(t_1, t_2) + 1.7 \cdot 10^5\frac{N}{\nu}+0.371 \cdot 10^5(t_2-t_1) \right. + \\ \nonumber && \qquad \qquad \qquad \left. +1+ |\omega|_{2+\frac{2}{3}\varepsilon}^{(1-\alpha)(2+\frac{2}{3}\varepsilon)}(t_1)\right)^{\frac{1}{1-\alpha}}
\end{eqnarray}
Therefore we obtain that
\begin{eqnarray}
&& \nu \int_{t_1}^{t_2} dt \int {\bar \nabla}(|{\bar \omega}({\bar x})|^{\frac{2}{3}\varepsilon}{\bar \omega}({\bar x}))({\bar \nabla}{\bar \omega}({\bar x}))d^3 x \leq \\ \nonumber && \leq 0.7305 \cdot 10^5 {\rm max}_{t \in [t_1, t_2]} \left(\int |{\bar \omega}({\bar x})|^{2+\frac{2}{3}\varepsilon} d^3x (t)\right)^{\alpha} \int_{t_1}^{t_2} dt \left( \int |{\bar \nabla}{\bar \omega}({\bar x})|^2 d^3x+1\right)  +\\ \nonumber && \qquad \qquad \qquad + \frac{1}{2+\frac{2}{3}\varepsilon}\left( \int |{\bar \omega}({\bar x})|^{2+\frac{2}{3}\varepsilon} d^3 x (t_1) + \int |{\bar \omega}({\bar x})|^{2+\frac{2}{3}\varepsilon} d^3 x (t_2)\right)+ \\ \nonumber && \qquad \qquad \qquad \qquad \qquad  + 2.1 {\rm max}_{t \in [t_1, t_2]} \int |{\bar \omega}({\bar x})|^{2+\frac{2}{3}\varepsilon} d^3 x N^{\frac{1}{2}} (t_2 - t_1)\\ \nonumber && 
\leq 1.97 e^{4.22 N^{\frac{1}{2}} (t_2-t_1)} \cdot \left(1+ \frac{3.99 \cdot 10^{26}}{\nu^{4+2\varepsilon}} N^{2+\varepsilon} C(t_1, t_2) + \frac{3.562 \cdot 10^5}{\nu} N^{\frac{1}{2}} C(t_1, t_2) + \right. \\ \nonumber && \left. \qquad \qquad \qquad+1.7 \cdot 10^5 \frac{N}{\nu}+0.371 \cdot10^5(t_2-t_1)  + 1.066 N^{\frac{1}{2}}(t_2-t_1) + |\omega|_{2+\frac{2}{3}\varepsilon}^{(1-\alpha)(2+\frac{2}{3}\varepsilon)}(t_1)\right)^5 
\end{eqnarray}
Therefore
\begin{eqnarray}
&& \frac{\frac{2}{3}\varepsilon}{(1+\frac{\varepsilon}{3})^2} \int_{t_1}^{t_2} dt \int({\bar \nabla}|{\bar \omega}({\bar x})|^{1+\frac{\varepsilon}{3}})^2 d^3 x + \int_{t_1}^{t_2} dt \int |{\bar \omega}({\bar x})|^{\frac{2}{3}\varepsilon} ({\bar \nabla}{\bar \omega}({\bar x}))^2 d^3 x \leq \qquad \qquad \qquad \\ \nonumber &&  \leq \frac{1.97}{\nu} e^{4.22 N^{\frac{1}{2}}(t_2-t_1)} \cdot \left(1+\frac{3.99 \cdot 10^{26}}{\nu^{4+2\varepsilon}} N^{2+\varepsilon} C(t_1, t_2) + \frac{3.562 \cdot 10^5}{\nu} N^{\frac{1}{2}} C(t_1, t_2)+\right. \\ \nonumber && \quad \quad \quad \left.+1.7 \cdot 10^5\frac{N}{\nu} + 0.371 \cdot 10^5(t_2-t_1) + 1.066 N^{\frac{1}{2}}(t_2-t_1) + |\omega|_{2+\frac{2}{3}\varepsilon}^{(1-\alpha)(2+\frac{2}{3}\varepsilon)}\right)^{5}. 
\end{eqnarray}

{\bf Part 3:} Differentiate with respect to $x$ the Navier-Stokes equation for ${\bar \omega}$.
\begin{eqnarray}
\partial_t \partial_x {\bar \omega}({\bar x}) - \nu \partial_x \Delta {\bar \omega} ({\bar x}) =  \partial_x \left(({\bar \omega}({\bar x}) \cdot {\bar \nabla}) {\bar u}({\bar x})\right) - \partial_x \left(({\bar u} ({\bar x}) \cdot {\bar \nabla}){\bar \omega}({\bar x})\right).
\end{eqnarray}

Take inner product with $\partial_x {\bar \omega}$
\begin{eqnarray}
\frac{1}{2} \partial_t(\partial_x {\bar \omega} ({\bar x}))^2 - \nu \partial_x {\bar \omega}({\bar x}) \cdot \partial_x \Delta {\bar \omega}({\bar x}) = \partial_x {\bar \omega}({\bar x}) \cdot \partial_x \left(({\bar \omega}({\bar x}) \cdot {\bar \nabla}){\bar u} ({\bar x})\right) -  \partial_x {\bar \omega}({\bar x}) \cdot \partial_x \left(({\bar u}({\bar x}) \cdot {\bar \nabla}) {\bar \omega}({\bar x})\right) \nonumber 
\end{eqnarray}

Integrate over space and by parts.
\begin{eqnarray}
&& \frac{1}{2} \partial_t \int (\partial_x {\bar \omega}({\bar x}))^2 d^3 x +  
\nu \int (\partial_{xx} {\bar \omega}({\bar x}))^2  d^3 x + \nu \int (\partial_{xy} {\bar \omega}({\bar x}))^2 d^3 x  + \nu \int (\partial_{xz} {\bar \omega}({\bar x}))^2 d^3 x = \qquad \nonumber \\ && \qquad \qquad \qquad  =- \int \partial_{xx} {\bar \omega}({\bar x}) \cdot ({\bar \omega}({\bar x}) \cdot {\bar \nabla}){\bar u}({\bar x}) d^3 x  + \int \partial_{xx} {\bar \omega} \cdot ({\bar u}({\bar x}) \cdot {\bar \nabla}){\bar \omega}({\bar x}) d^3 x.
\end{eqnarray}
Add the equations obtained similarly but with differentiation with respect to $y$ or $z$ instead of $x$. 
\begin{eqnarray}
&& \frac{1}{2} \partial_t \int (\partial_y {\bar \omega}({\bar x}))^2 d^3 x + \nu \int (\partial_{yy} {\bar \omega}({\bar x}))^2 d^3 x + \nu \int (\partial_{xy} {\bar \omega}({\bar x}))^2 d^3 x + \nu \int (\partial_{yz} {\bar \omega}({\bar x}))^2 d^3 x  = \qquad \nonumber \\ && \qquad \qquad \qquad   =- \int \partial_{yy} {\bar \omega}({\bar x}) \cdot ({\bar \omega}({\bar x}) \cdot {\bar \nabla}) {\bar u}({\bar x}) d^3 x + \int \partial_{yy} {\bar \omega}({\bar x}) \cdot ({\bar u}({\bar x}) \cdot {\bar \nabla}) {\bar \omega}({\bar x}) d^3 x.
\end{eqnarray}
and 
\begin{eqnarray}
&&\frac{1}{2} \partial_t \int (\partial_z {\bar \omega}({\bar x}))^2 d^3 x +\nu \int (\partial_{zz} {\bar \omega}({\bar x}))^2 d^3 x + \nu \int (\partial_{xz} {\bar \omega}({\bar x}))^2 d^3 x + \nu \int (\partial_{yz} {\bar \omega}({\bar x}))^2 d^3 x = \qquad \nonumber \\ &&  \qquad \qquad \qquad  = -\int \partial_{zz} {\bar \omega}({\bar x}) \cdot ({\bar \omega}({\bar x}) \cdot {\bar \nabla}) {\bar u}({\bar x}) d^3 x + \int \partial_{zz} {\bar \omega}({\bar x}) \cdot ({\bar u}({\bar x}) \cdot {\bar \nabla}) {\bar \omega} ({\bar x}) d^3 x.
\end{eqnarray}
to get
\begin{eqnarray}
&& \frac{1}{2} \partial_t \int |\nabla \omega|^2 ({\bar x}) d^3 x + \nu \int |{\bar \nabla} {\bar \nabla} \omega|^2 ({\bar x}) d^3 x =\\ \nonumber && \qquad \qquad \qquad = - \int \Delta {\bar \omega}({\bar x}) \cdot ({\bar \omega}({\bar x}) \cdot {\bar \nabla}) {\bar u}({\bar x}) d^3 x + \int \Delta {\bar \omega}({\bar x}) \cdot ({\bar u}({\bar x}) \cdot {\bar \nabla}) {\bar \omega}({\bar x}) d^3 x.  \nonumber 
\end{eqnarray}
Plugging the expression for ${\bar u}$ in terms of ${\bar \omega}$ 
$$ {\bar u}({\bar x}) = \frac{1}{4 \pi} \int \frac{{\bar \omega}({\bar y}) \times {({\bar x}-{\bar y})}}{|x-y|^3} d^3 y = \frac{1}{4\pi} \int \frac{{\bar \omega}({\bar x} + {\bar y}) \times {\hat y}}{|y|^2} d^3 y$$ on the right hand side, we obtain
\begin{eqnarray}\label{eqn35}
&& \frac{1}{2} \partial_t \int |{\bar \nabla} {\bar \omega}({\bar x})|^2 d^3 x + \nu \int |{\bar \nabla}{\bar \nabla} {\bar \omega}({\bar x})|^2 d^3 x = - \frac{1}{4 \pi} \int \int \Delta {\bar \omega}({\bar x}) \cdot ({\bar \omega}({\bar x}) \cdot {\bar \nabla}) \frac{{\bar \omega}({\bar x}+{\bar y}) \times {\hat y}}{|y|^2} d^3 x d^3 y \nonumber \\ && \qquad \qquad \qquad \qquad \qquad  + \frac{1}{4 \pi} \int \int \Delta {\bar \omega}({\bar x}) \cdot \frac{(({\bar \omega}({\bar x} + {\bar y}) \times {\hat y}) \cdot {\bar \nabla}) {\bar \omega}({\bar x})}{|y|^2} d^3x d^3 y. 
\end{eqnarray}
We can rewrite the right hand side of ({\ref{eqn35}}) as
\begin{eqnarray}
&& -\frac{1}{4\pi} \int \int_{|y| \leq 1} \Delta {\bar \omega}({\bar x}) \cdot ({\bar \omega} ({\bar x}) \cdot {\bar \nabla})\frac{{\bar \omega}({\bar x} + {\bar y}) \times {\hat y}}{|y|^2} d^3 x d^3 y-\\ \nonumber && \qquad \qquad \qquad - \frac{1}{4 \pi} \int \int_{|y| > 1}
\Delta {\bar \omega}({\bar x}) \cdot ({\bar \omega}({\bar x}) \cdot {\bar \nabla}) \frac {{\bar \omega}({\bar x} + {\bar y}) \times {\hat y}}{|y|^2} d^3 x d^3 y +\\ \nonumber && + \frac{1}{4 \pi} \int \int_{|y| \leq 1} \Delta {\bar \omega}({\bar x}) \cdot \frac{(({\bar \omega}({\bar x} + {\bar y}) \times {\hat y}) \cdot {\bar \nabla}) {\bar \omega}({\bar x})}{|y|^2} d^3 x d^3 y+ \\ \nonumber && \qquad \qquad \qquad + \frac{1}{4 \pi} \int \int_{|y| > 1} \Delta {\bar \omega}({\bar x}) \cdot \frac{(({\bar \omega}({\bar x} + {\bar y}) \times {\hat y}) \cdot {\bar \nabla}) {\bar \omega}({\bar x})}{|y|^2} d^3 x d^3 y.
\end{eqnarray}
The absolute value of the integrals over $|y| \leq 1$ can be estimated from above by
\begin{eqnarray}\nonumber
&& \frac{3}{2 \pi} \int \int_{|y| \leq 1} |\Delta {\bar \omega}({\bar x})| |{\bar \omega}({\bar x})| \frac{|{\bar \nabla}{\bar \omega}({\bar x}+{\bar y})|}{|y|^2} d^3x d^3 y +  
\frac{3}{2 \pi} \int \int_{|y| \leq 1} |\Delta {\bar \omega}({\bar x})| |{\bar \omega}({\bar x} + {\bar y})| \frac{|{\bar \nabla} {\bar \omega}({\bar x})|}{|y|^2} d^3 x d^3 y \\ \nonumber && = \frac{3}{2 \pi} \int \int_{|y| \leq 1} |\Delta {\bar \omega}({\bar x} - {\bar y})| |{\bar \omega} ({\bar x} - {\bar y})| \frac{|{\bar \nabla} {\bar \omega}({\bar x})|}{|y|^2} d^3 xd^3 y \\ \nonumber && \qquad \qquad \qquad \qquad \qquad + \frac{3}{2 \pi} \int \int_{|y| \leq 1} |\Delta {\bar \omega}({\bar x})| |{\bar \omega}({\bar x} + {\bar y})| \frac{|{\bar \nabla}{\bar \omega}({\bar x})|}{|y|^2} d^3 x d^3 y \\ \nonumber && =  \frac{3}{2\pi} \int_{|y| \leq 1} \frac{d^3 y}{|y|^2} \int |\Delta {\bar \omega}({\bar x}-{\bar y})| |{\bar \omega}({\bar x}-{\bar y})| |{\bar \nabla}{\bar \omega}({\bar x})| d^3 x + \\ \nonumber && \qquad \qquad \qquad  + \frac{3}{2 \pi} \int |\Delta {\bar \omega}({\bar x})| |{\bar \nabla}{\bar \omega}({\bar x})| \int_{|y| \leq 1} \frac{|{\bar \omega}({\bar x}+{\bar y})|}{|y|^2} d^3 y \\ \nonumber && \leq \frac{3}{4 \pi} \int_{|y| \leq 1} \frac{d^3 y}{|y|^2} \int (\frac{\nu}{5} |\Delta {\bar \omega}({\bar x} - {\bar y})|^2 + \frac{5}{\nu}|{\bar \omega}({\bar x}-{\bar y})|^2 |{\bar \nabla}{\bar \omega}({\bar x})|^2) d^3 x + \\ \nonumber &&  + \frac{3}{4 \pi} \int(\frac{\nu}{5} |\Delta {\bar \omega}({\bar x})|^2 + \frac{5}{\nu} |{\bar \nabla} {\bar \omega}({\bar x})|^2 \int_{|y| \leq 1} \frac{|{\bar \omega}({\bar x}+ {\bar y})|}{|y|^2} d^3y \int_{|z| \leq 1} \frac{|{\bar \omega}({\bar x} + {\bar z})|}{|z|^2} d^3 z) d^3 x\\ \nonumber && = \frac{3}{4 \pi} \frac{\nu}{5} \int |\Delta {\bar \omega}({\bar x})|^2 (1 + \int_{|y| \leq 1} \frac{d^3 y}{|y|^2}) + \frac{3}{4 \pi} \frac{5}{\nu} \int |{\bar \nabla} {\bar \omega}({\bar x})|^2 d^3 x \int_{|y| \leq 1} \frac{|{\bar \omega}({\bar x}+{\bar y})|^2}{|y|^2} d^3 y  (\int \frac{d^3 z}{|z|^2} + 1) \\ \nonumber && \leq \frac{3}{4 \pi}\frac{\nu}{5} \int|\Delta {\bar \omega}({\bar x})|^2 d^3 x(1+ \int_{|y| \leq 1} \frac{d^3 y}{|y|^2}) + \\ \nonumber && \qquad \qquad \qquad + \frac{3}{4 \pi} \frac{5}{\nu} \int |{\bar \nabla} {\bar \omega}({\bar x})|^2 d^3 x (\int |\omega|^{6+2\varepsilon}({\bar x}) d^3 x)^{\frac{1}{3+\varepsilon}}(4 \pi + 1) (4 \pi)^{\frac{2+\varepsilon}{3+\varepsilon}}(\frac{2+\varepsilon}{\varepsilon})^{\frac{2+\varepsilon}{3+\varepsilon}} \\ \nonumber && \leq \frac{3(1+4\pi)}{4 \pi}\frac{\nu}{5} \int |\Delta {\omega}({\bar x})|^2 d^3 x+ \frac{30}{\nu \varepsilon} (4 \pi +1)\int|{\bar \nabla}{\bar \omega}({\bar x})|^2 d^3 x (\int |{\bar \nabla}|{\bar \omega}|^{1+\frac{\varepsilon}{3}}({\bar x})|^2 d^3 x)^{\frac{3}{3+\varepsilon}}.
\end{eqnarray}   
The absolute value of the integral over $|y| > 1$ is bounded by 
\begin{eqnarray}
&& \frac{1}{4 \pi} \int \int_{|y| > 1} \left|\Delta {\bar \omega}({\bar x}) \cdot ({\bar \omega}({\bar x}) \cdot {\bar \nabla}) \frac{{\bar \omega}({\bar x}+{\bar y}) \times {\hat y}}{|y|^2}\right| d^3 xd^3 y + \\ \nonumber && \qquad \qquad \qquad \qquad \qquad + \frac{1}{4 \pi}\int \int_{|y| > 1} \left|\Delta {\bar \omega}({\bar x}) \cdot \frac{(({\bar \omega}({\bar x}+{\bar y}) \times {\hat y}) \cdot {\bar \nabla}) {\bar \omega}({\bar x})}{|y|^2}\right| d^3 x d^3 y \\ \nonumber && \leq
\frac{3}{2 \pi} \int \int_{|y| > 1} |\Delta {\bar \omega}({\bar x})| |{\bar \omega}({\bar x})| |{\bar \nabla}{\bar \omega}({\bar x} + {\bar y})| d^3 x \frac{d^3 y}{|y|^2} +\\ \nonumber && \qquad \qquad \qquad \qquad \qquad 
+ \frac{3}{2\pi} \int |\Delta {\bar \omega}({\bar x})| |{\bar \nabla}{\bar \omega}({\bar x})| \int_{|y| \geq 1} \frac{|{\bar \omega}({\bar x}+{\bar y})|}{|y|^2} d^3 y d^3x
\\ \nonumber && \leq \frac{3}{4 \pi} \frac{\nu}{5} \int |\Delta {\bar \omega}({\bar x})|^2 d^3 x  +  \frac{15}{\nu} \int |\omega|^2({\bar x}) d^3 x (\int |{\bar \nabla}{\bar \omega}|^2({\bar x} + {\bar y}) d^3 y) +\\ \nonumber && \qquad \qquad \qquad + 
\frac{3}{4\pi} \frac{\nu}{5} \int |\Delta{\bar \omega}({\bar x})|^2 d^3 x +  \frac{15}{\nu} \int |{\bar \nabla} {\bar \omega}({\bar x})|^2 (\int |{\bar \omega}({\bar x}+{\bar y})|^2 d^3 y) d^3x.\end{eqnarray}
Therefore we obtain from ({\ref{eqn35}}) that 
\begin{eqnarray}
&& \frac{1}{2} \partial_t \int |{\bar \nabla} {\bar \omega}({\bar x})|^2 d^3 x + \frac{(8 \pi - 9) \nu}{20 \pi} \int |{\bar \nabla}{\bar \nabla} {\bar \omega}({\bar x})|^2 d^3x \leq \\ \nonumber &&  \leq \frac{30(4\pi+1)}{\varepsilon \nu} \int |{\bar \nabla}{\bar \omega}({\bar x})|^2 d^3x (\int ({\bar \nabla}|{\bar \omega}({\bar x})|^{1+\frac{\varepsilon}{3}})^2 d^3 x)^{\frac{3}{3+\varepsilon}}+ \frac{30}{\nu} \int |{\bar \nabla}{\bar \omega}({\bar x})|^2 d^3 x \int |{\bar \omega}({\bar x})|^2 d^3 x. \nonumber
\end{eqnarray}
We derive that
\begin{eqnarray}
 \partial_t {\rm ln} \int |{\bar \nabla}{\bar \omega}({\bar x})|^2 d^3 x \leq \frac{816}{\varepsilon \nu} (\int |{\bar \nabla}|{\bar \omega}({\bar x})|^{1+\frac{\varepsilon}{3}}|^2 d^3 x+1)+\frac{60}{\nu} \int |{\bar \omega}({\bar x})|^2 d^3 x.\qquad
\end{eqnarray}
Therefore
\begin{eqnarray}
&& {\rm ln} \int |{\bar \nabla}{\bar \omega}({\bar x})|^2 d^3 x(t_2) - {\rm ln} \int |{\bar \nabla}{\bar \omega}({\bar x})|^2 d^3 x(t_1) \leq \nonumber  \\ \nonumber && \leq \frac{8.16 \cdot 10^4}{\nu} \left(\int_{t_1}^{t_2} dt \int |{\bar \nabla}|{\bar \omega}({\bar x})|^{1+\frac{\varepsilon}{3}}|^2 d^3 x (t)+(t_2-t_1)\right) + \frac{60}{\nu} C(t_1, t_2)\ \ \   \nonumber \\ && \leq  \frac{8.16 \cdot 10^4}{\nu} (t_2-t_1) + \frac{60}{\nu} C(t_1, t_2) + \frac{2.432 \cdot 10^{7}}{\nu^2} e^{4.22 N^{\frac{1}{2}}(t_2-t_1)} \times \\ \nonumber && \left(1+ \frac{3.99 \cdot 10^{26}}{\nu^{4.02}} N^{2.01} C(t_1, t_2) + \frac{3.562 \cdot 10^5}{\nu} N^{\frac{1}{2}} C(t_1, t_2) + 1.7 \cdot 10^5 \frac{N}{\nu} +0.371 \cdot10^5(t_2-t_1)\right. \\ \nonumber && \qquad \qquad \qquad \qquad \qquad \left. + 1.066 N^{\frac{1}{2}}(t_2-t_1) + |{\bar \omega}|^{(2+\frac{2}{3} \varepsilon)(1-\alpha)}_{2+\frac{2}{3}\varepsilon}\right)^{5} 
\end{eqnarray}
This implies
\begin{eqnarray}
&& \int |{\bar \nabla}{\bar \omega}({\bar x})|^2 d^3x (t_2) \leq  \int |{\bar \nabla}{\bar \omega}({\bar x})|^2 d^3 x (t_1) {\rm exp}\left(\frac{8.16 \cdot 10^4}{\nu} (t_2-t_1) + \frac{60}{\nu} C(t_1,t_2)+  \right. \qquad \qquad \qquad  \\ \nonumber && \left. + \frac{2.432 \cdot 10^{7}}{\nu^2}e^{4.22 N^{\frac{1}{2}} (t_2-t_1)}\left(1 + \frac{3.99 \cdot 10^{26}}{\nu^{4.02}} N^{2.01} C(t_1, t_2)+\frac{3.562 \cdot 10^5}{\nu} N^{\frac{1}{2}} C(t_1, t_2) + \right. \right. \\ \nonumber && \left. \left.  + 1.7 \cdot 10^{5}\frac{N}{\nu} + 0.371 \cdot 10^5 (t_2-t_1) + 1.066 N^{\frac{1}{2}}(t_2-t_1) + |\omega|_{2+\frac{2}{3}\varepsilon}^{(2+\frac{2}{3}\varepsilon)(1-\alpha)}(t_1)\right)^{5}\right)\nonumber
\end{eqnarray} 
We obtain 
\begin{eqnarray}
&& \nonumber \int_{t_1}^{t_2} dt \int |{\bar \nabla} {\bar \nabla} {\bar \omega}({\bar x})|^2 d^3 x \leq \frac{8.16 \ \pi \ 10^5}{(8\pi-9) \nu^2} \int_{t_1}^{t_2} dt (\int |{\bar \nabla}{\bar \omega}({\bar x})|^2 d^3 x) (\int ({\bar \nabla}|{\bar \omega}({\bar x})|^{1+\frac{\varepsilon}{3}})^2 d^3x)^{\frac{3}{3+\varepsilon}}  \\ \nonumber && \quad \quad \quad  + \frac{10 \pi}{(8 \pi -9) \nu} \int |{\bar \nabla}{\bar \omega}({\bar x})|^2 d^3 x (t_1) + \frac{10 \pi}{(8 \pi - 9) \nu} \int |{\bar \nabla}{\bar \omega}({\bar x})|^2 d^3 x (t_1)+ \\ \nonumber && \quad \quad \quad \quad \quad + \frac{6 \cdot 10^2 \pi}{(8\pi-9) \nu} \int_{t_1}^{t_2} dt \int |{\bar \nabla}{\bar \omega}({\bar x})|^2 d^3 x \int |{\bar \omega}({\bar x})|^2 d^3 x \\ \nonumber &&
\leq \frac{1.945}{\nu} \int |{\bar \nabla}{\bar \omega}({\bar x})|^2 d^3 x (t_1)\ 
{\rm exp}\left(\frac{1.632 \cdot 10^5(t_2 -t_1)}{\nu} + 120 \frac{C(t_1, t_2)}{\nu}+ \right. \\ \nonumber &&\left.+ \frac{4.864 \cdot 10^{7}}{\nu^2}e^{4.22 N^{\frac{1}{2}} (t_2-t_1)}\cdot \left(1 + \frac{3.99 \cdot 10^{26}}{ \nu^{4.02}} N^{2+\varepsilon} C(t_1, t_2)+\frac{3.562 \cdot 10^5}{\nu}N^{\frac{1}{2}} C(t_1,t_2)+ \right.\right. \\ \nonumber && \left. \left.\qquad \qquad \qquad + 1.7 \cdot 10^5\frac{N}{\nu}+0.371 \cdot 10^5(t_2-t_1) + 1.066 N^{\frac{1}{2}}(t_2-t_1)+ |\omega|^{(2+\frac{2}{3}\varepsilon)(1-\alpha)}_{2+\frac{2}{3}\varepsilon}(t_1) \right)^{5}\right).
\end{eqnarray}

{\bf Part 4:} Now we are going to estimate $\int |{\bar \nabla}{\bar \omega}|^{3+\varepsilon} ({\bar x}) d^3 x $.

Taking the gradient of the Navier-Stokes equation for $\omega$ and then taking the scalar product with ${\bar \nabla} {\bar \omega}$ we obtain
\begin{eqnarray}\label{ns29}
&& \frac{1}{2} \partial_t ({\bar \nabla} {\bar \omega})^2 - \nu \sum_i \partial_i {\bar \omega}({\bar x}) \cdot \partial_i (\Delta {\bar \omega}({\bar x})) = \sum_i \partial_i {\bar \omega}({\bar x}) \cdot \partial_i \left(({\bar \omega}({\bar x}) \cdot {\bar \nabla}){\bar u}({\bar x})\right)- \qquad \qquad \qquad \\ \nonumber &&\qquad \qquad \qquad \qquad \qquad  -  \sum_i \partial_i {\bar \omega}({\bar x}) \cdot \partial_i \left(({\bar u}({\bar x}) \cdot {\bar \nabla}){\bar \omega}({\bar x})\right).
\end{eqnarray}
Multiplying ({\ref{ns29}) by $|{\bar \nabla} {\bar \omega}|^{1+\varepsilon}$ and integrating we obtain
\begin{eqnarray}\label{ns30}
&& \frac{1}{3+\varepsilon} \partial_t \int |{\bar \nabla} {\bar \omega}({\bar x})|^{3+\varepsilon} d^3 x - \nu \int |{\bar \nabla} {\bar \omega}|^{1+\varepsilon}({\bar x}) \sum_i \partial_i {\bar \omega}({\bar x}) \cdot \partial_i (\Delta {\bar \omega}({\bar x})) d^3x \\ \nonumber && \qquad \qquad \qquad  = \int |{\bar \nabla}{\bar \omega}({\bar x})|^{1+\varepsilon} \sum_i \partial_i {\bar \omega}({\bar x}) \cdot \partial_i \left(({\bar \omega}({\bar x})\cdot {\bar \nabla})\int \frac{{\bar \omega}({\bar x} +{\bar y}) \times {\hat y}}{|y|^2}\right) d^3 y d^3x - \\ \nonumber && \qquad \qquad \qquad \qquad \qquad - \int |{\bar \nabla}{\bar \omega}({\bar x})|^{1+\varepsilon} \sum_i \partial_i {\bar \omega}({\bar x}) \cdot \partial_i \left((\int \frac{{\bar \omega}({\bar x}+{\bar y}) \times {\hat y}}{|y|^2} d^3 y \cdot {\bar \nabla}) {\bar \omega}({\bar x})\right) d^3x. 
\end{eqnarray}
The viscosity term in (\ref{ns30}) can be rewritten as
\begin{eqnarray}
&& - \nu \int d^3 x |{\bar \nabla}{\bar \omega}({\bar x})|^{1+\varepsilon} \sum_i \partial_i {\bar \omega}({\bar x}) \cdot \partial_i (\sum_j \partial_j \partial_j {\bar \omega}({\bar x}))  = \\ \nonumber && 
= \nu \int d^3 x |{\bar \nabla} {\bar \omega}({\bar x})|^{1+\varepsilon} \sum_i \sum_j (\partial_i \partial_j {\bar \omega} ({\bar x}))^2 + \nu \int d^3 x \frac{1+\varepsilon}{|{\bar \nabla}{\bar \omega}({\bar x})|^{1-\varepsilon}} \sum_j (\sum_i \partial_i {\bar \omega}({\bar x}) \cdot \partial_i \partial_j {\bar \omega}({\bar x}))^2.
\end{eqnarray}
Therefore the viscosity term in the equation (\ref{ns30}) is positive definite.

 The absolute value of the right hand side of ({\ref{ns30}}) can be rewritten as
\begin{eqnarray}\label{ns31}
&& \left|\sum_{i,j,k} \int \int d^3 x  d^3 y |{\bar \nabla}{\bar \omega}({\bar x})|^{1+\varepsilon} \partial_i \omega_k ({\bar x}) \partial_i \omega_j ({\bar x}) \ \frac{\partial_j ({\bar \omega}({\bar x}+{\bar y}) \times {\hat y})_k}{|y|^2} \right. + \\ \nonumber && \quad \quad \quad  + \sum_{i,j,k} \int \int d^3 x d^3 y |{\bar \nabla} {\bar \omega}({\bar x})|^{1+\varepsilon} \partial_i \omega_k ({\bar x}) \omega_j({\bar x}) \frac{ \partial_i \partial_j ({\bar \omega}({\bar x}+{\bar y}) \times {\hat y})_k}{|y|^2} -\\ \nonumber &&  \qquad \qquad \qquad - \left. \sum_{i,j,k} \int \int d^3 x d^3 y 
|{\bar \nabla} {\bar \omega}({\bar x})|^{1+\varepsilon} \partial_i \omega_k ({\bar x}) \frac{ \partial_i ({\bar \omega}({\bar x}+{\bar y}) \times {\hat y})_j}{|y|^2} \ \partial_j \omega_k({\bar x})\right. - \\ \nonumber && \qquad \qquad \qquad \qquad \qquad - \left. \sum_{i,j,k} \int \int d^3 x d^3 y | {\bar \nabla}{\bar \omega}({\bar x})|^{1+\varepsilon} \partial_i \omega_k({\bar x})  \frac{({\bar \omega}({\bar x}+{\bar y}) \times {\hat y})_j}{|y|^2}\ \partial_i \partial_j \omega_k({\bar x})\right| \leq \\ \nonumber && \leq
12 \sqrt{3} \int \int d^3 x d^3 y |{\bar \nabla} {\bar \omega}({\bar x})|^{3+\varepsilon} \frac{|{\bar \nabla}{\bar \omega}({\bar x}+{\bar y})|}{|y|^2} + \\ \nonumber && \quad \quad \quad \qquad \qquad \qquad  +6 \sqrt{3}\int d^3 x \int_{|y| \leq 1}
d^3 y |{\bar \nabla}{\bar \omega}({\bar x})|^{2+\varepsilon} |{\bar \omega}({\bar x})| \frac{|{\bar \nabla}{\bar \nabla}{\bar \omega}({\bar x}+{\bar y})|}{|y|^2} + \\ \nonumber && \quad \quad \quad + 6 \sqrt{3} \int d^3 x \int_{|y| > 1} d^3 y |{\bar \nabla}{\bar \omega}({\bar x})|^{2+\varepsilon} |{\bar \omega}({\bar x})| \frac{|{\bar \nabla}{\bar \nabla}{\bar \omega}({\bar x}+{\bar y})|}{|y|^2}+ \\ \nonumber &&\quad \quad \quad  \qquad \qquad \qquad + 6 \sqrt{3} \int \int_{|y| \leq 1} d^3 x d^3 y |{\bar \nabla}{\bar \omega}|^{2+\varepsilon}({\bar x}) |{\bar \omega}({\bar x}+{\bar y})|\frac{|{\bar \nabla}{\bar \nabla}{\bar \omega}({\bar x})|}{|y|^2} + \\ \nonumber && \quad \quad \quad \qquad \qquad \qquad \qquad \qquad +6 \sqrt{3}\int \int_{|y| > 1} d^3 x d^3 y |{\bar \nabla}{\bar \omega}({\bar x})|^{2+\varepsilon} |{\bar \omega}({\bar x}+{\bar y})| \frac{|{\bar \nabla}{\bar \nabla}{\bar \omega}({\bar x})|}{|y|^2}. 
\end{eqnarray} 

The first integral on the right hand side of ({\ref{ns31}}) can be estimated as
\begin{eqnarray}
&& \int \int_{|y| \leq 1}  d^3 x d^3 y |{\bar \nabla}{\bar \omega}({\bar x})|^{3+\varepsilon} \frac{|{\bar \nabla}{\bar \omega}({\bar x}+{\bar y})|}{|y|^2} + \int \int_{|y| > 1} d^3 x d^3y |{\bar \nabla}{\bar \omega}({\bar x})|^{3+\varepsilon} \frac{|{\bar \nabla}{\bar \omega}({\bar x}+{\bar y})|}{|y|^2} \\ \nonumber &&  \leq \int d^3 x |{\bar \nabla}{\bar \omega}({\bar x})|^{3+\varepsilon}(\int |{\bar \nabla}{\bar \omega}({\bar x}+{\bar y})|^{6} d^3 y)^{\frac{1}{6}}(\frac{20 \pi}{3})^{\frac{5}{6}} + \int d^3 x |{\bar \nabla}{\bar \omega}({\bar x})|^{3+\varepsilon} \int d^3 y \frac{|{\bar \nabla}{\bar \omega}({\bar x}+{\bar y})|}{|y|^2} \\ \nonumber && \leq \frac{20 \pi}{3}\int d^3 x |{\bar \nabla}{\bar \omega}({\bar x})|^{3+\varepsilon} (\int d^3 y |{\bar \nabla} |{\bar \nabla} {\bar \omega}({\bar x}+{\bar y})||^{2})^{\frac{1}{2}} + \int d^3 x |{\bar \nabla}{\bar \omega}({\bar x})|^{3+\varepsilon} \int d^3y \frac{ |{\bar \nabla}{\bar \omega}({\bar x}+{\bar y})|}{|y|^2}\leq  \\ \nonumber && \leq 
\frac{20 \pi}{\sqrt{3}} \int d^3 x |{\bar \nabla}{\bar \omega}({\bar x})|^{3+\varepsilon} (\int d^3y |{\bar \nabla}{\bar \nabla}{\bar \omega}({\bar y})|^2)^{\frac{1}{2}} + \sqrt{4 \pi} \int d^3 x |{\bar \nabla}{\bar \omega}({\bar x})|^{3+\varepsilon} (\int d^3 y |{\bar \nabla}{\bar \omega}({\bar x}+{\bar y})|^2)^{\frac{1}{2}}\leq \\ \nonumber && \leq \frac{20 \pi}{\sqrt {3}} \int d^3 x |{\bar \nabla}{\bar \omega}|^{3+\varepsilon} (\int d^3 y |{\bar \nabla}{\bar \nabla}{\bar \omega}|^2 +1) + \sqrt{4 \pi} \int d^3 x |{\bar \nabla}{\bar \omega}({\bar x})|^{3+\varepsilon} (\int d^3 y |{\bar \nabla}{\bar \omega}({\bar y})|^2 +1).
\end{eqnarray}

The second integral on the right hand side of (\ref{ns31}) can be evaluated as
\begin{eqnarray}
&& \int d^3 x \int_{|y| \leq 1} d^3 y |{\bar \nabla} {\bar \omega}({\bar x})|^{2+\varepsilon} |{\bar \omega}({\bar x})| \frac{|{\bar \nabla}{\bar \nabla}{\bar \omega}({\bar x}+{\bar y})|}{|y|^2}\\ \nonumber &&  = \int d^3 x \int_{|y| \leq 1} d^3 y |{\bar \nabla} {\bar \omega} ({\bar x} - {\bar y})|^{2+\varepsilon} |{\bar \omega}({\bar x}-{\bar y})| \frac{|{\bar \nabla}{\bar \nabla}{\bar \omega}({\bar x})|}{|y|^2}  \leq \sqrt{\int |{\bar \nabla}{\bar \nabla} {\bar \omega}({\bar x})|^2 d^3 x} \times \\ \nonumber &&  \times \sqrt{ \int d^3 x \int_{|y| \leq 1} d^3 y \frac{|{\bar \nabla}{\bar \omega}({\bar x}-{\bar y})|^{2+\varepsilon} |{\bar \omega}({\bar x} - {\bar y})|}{|y|^2} \int_{|z| \leq 1} d^3 z \frac{|{\bar \nabla}{\bar \omega}({\bar x}-{\bar z})|^{2+\varepsilon} |{\bar \omega}({\bar x} - {\bar z})|}{|z|^2}} \\ \nonumber && \leq \sqrt{\int |{\bar \nabla}{\bar \nabla} {\bar \omega}({\bar x})|^2 d^3 x }\times \\ \nonumber && \times  \sqrt{\int d^3 x \int_{|y| \leq 1} d^3 y \frac{|{\bar \nabla}{\bar \omega}({\bar x} - {\bar y})|^{2+\varepsilon}|{\bar \omega}({\bar x}-{\bar y})|}{|y|^{3-\frac{\varepsilon}{2}}} \int_{|z| \leq 1} d^3 z \frac{|{\bar \nabla}{\bar \omega}({\bar x} - {\bar z})|^{2+\varepsilon}|{\bar \omega}({\bar x}-{\bar z})|}{|z|^{1+\frac{\varepsilon}{2}}}} \\ \nonumber && \leq \sqrt{\int |{\bar \nabla} {\bar \nabla} {\bar \omega}({\bar x})|^2 d^3 x} (\int \frac{d^3 z}{|z|^{3(1-\frac{\varepsilon^2}{4})}})^{\frac{1}{6(1-\frac{\varepsilon}{2})}} \times \\ \nonumber && \sqrt{\int d^3 x' |{\bar \nabla}{\bar \omega}|^{2+\varepsilon} |{\bar \omega}({\bar x'})| \int_{|y| \leq 1} \frac{d^3 y}{|y|^{3-\frac{\varepsilon}{2}}} \left(\int d^3 z  |{\bar \nabla}{\bar \omega}({\bar x'} + {\bar y} - {\bar z}) |^{\frac{3(1-\frac{\varepsilon^2}{4})}{1-\frac{3\varepsilon}{4}}} |{\bar \omega}({\bar x'}+{\bar y} - {\bar z})|^{\frac{3(1-\frac{\varepsilon}{2})}{2(1-\frac{3}{4}\varepsilon)}}\right)^{\frac{2-\frac{3}{2}\varepsilon}{3 - \frac{3}{2}\varepsilon}}}\\ \nonumber && \leq  (\frac{16 \pi}{3 \varepsilon^2})^{\frac{1}{6 - 3 \varepsilon}} (\frac{8 \pi}{\varepsilon})^{\frac{1}{2}} \sqrt{\int |{\bar \nabla}{\bar \nabla} {\bar \omega}({\bar x})|^2 d^3 x} \sqrt{\int d^3 x |{\bar \omega}({\bar x})|^{2+\varepsilon} |{\bar \omega}({\bar x})|} \\ \nonumber && \qquad \qquad \qquad \cdot \left(\int d^3 z |{\bar \nabla}{\bar \omega}|^{\frac{3(1-\frac{\varepsilon^2}{4})}{1-\frac{3}{4}\varepsilon}} |{\bar \omega}({\bar z})|^{\frac{3(1-\frac{\varepsilon}{2})}{2(1-\frac{3}{4}\varepsilon)}}\right)^{\frac{1}{3} \frac{1-\frac{3}{4}\varepsilon}{1-\frac{\varepsilon}{2}}}  \\ \nonumber &&  \leq 629.2 (\int |{\bar \nabla}{\bar \nabla}{\bar  \omega}({\bar x})|^2 d^3 x)^{\frac{1}{2}} (\int |{\bar \nabla}{\bar \omega}({\bar x})|^{3+\varepsilon} d^3 x)^{\frac{1}{3}\frac{1+\frac{\varepsilon}{2}}{1+\frac{\varepsilon}{3}}} (\int |{\bar \nabla}{\bar \omega}({\bar x})|^2 d^3 x)^{\frac{1}{8}\frac{1+\varepsilon}{1+\frac{\varepsilon}{3}}} (\int |\omega|^2({\bar x}) d^3 x)^{\frac{1}{8} \frac{1-\frac{\varepsilon}{3}}{1+\frac{\varepsilon}{3}}} \\ \nonumber &&
(\int |{\bar \nabla}{\bar \nabla}{\bar \omega}({\bar x})|^2 d^3 x)^{\frac{1}{4} \frac{1+\frac{3}{2}\varepsilon -\frac{\varepsilon^2}{6}}{1-\frac{7}{6}\varepsilon+\frac{5 \varepsilon^2}{18}} \frac{1-\frac{5}{6}\varepsilon}{1-\frac{\varepsilon}{2}}} (\int |{\bar \nabla}{\bar \omega}({\bar x})|^{3+\varepsilon} d^3x)^{\frac{1}{6} \frac{1 -\frac{5}{2}\varepsilon + \frac{\varepsilon^2}{2}}{1-\frac{7}{6}\varepsilon +\frac{5 \varepsilon^2}{18}} \frac{1-\frac{5}{6}\varepsilon}{1 -\frac{\varepsilon}{2}}} (\int |{\bar \nabla}{\bar \omega}|^2({\bar x}) d^3 x)^{\frac{1}{4}} \\ \nonumber && \leq 1.1 \cdot 10^3 (\int |{\bar \nabla}{\bar \nabla}{\bar \omega}({\bar x})|^2 d^3 x +1) (\int |{\bar \nabla}{\bar \omega}({\bar x})|^{3+\varepsilon} d^3x +1) (\int |{\bar \nabla}{\bar \omega}({\bar x})|^2 d^3 x + 1) (N+1).
\end{eqnarray}
The third integral on the right hand side of the equation (\ref{ns30}) can be evaluated as
\begin{eqnarray}
&& \int \int_{|y| > 1} |{\bar \nabla}{\bar \omega}({\bar x})|^{2+\varepsilon} |{\bar \omega}({\bar x})| \frac{|{\bar \nabla}{\bar \nabla} {\bar \omega}({\bar x}+{\bar y})|}{|y|^2} d^3 x d^3 y \\ \nonumber && \leq \int d^3x  |{\bar \nabla}{\bar \omega}({\bar x})|^{3+\varepsilon} |{\bar \omega}({\bar x})| \sqrt{\int d^3 y |{\bar \nabla}{\bar \nabla}{\bar \omega}({\bar x} + {\bar y})|^2 d^3 y} \sqrt{\int_{|y| > 1}\frac{d^3 y}{|y|^4}} \\ \nonumber && \leq \sqrt{4 \pi} \sqrt{\int d^3 y |{\bar \nabla}{\bar \nabla} {\bar \omega}({\bar y})|^2} \int d^3x |{\bar \nabla}{\bar \omega}({\bar x})|^{2+\varepsilon} |{\bar \omega}({\bar x})| \\ \nonumber && \leq \sqrt{4 \pi} \sqrt{\int d^3 x |{\bar \nabla}{\bar \nabla}{\bar \omega}({\bar x})|^2} (\int d^3 x  |{\bar \nabla} {\bar \omega}({\bar x})|^{3+\varepsilon})^{\frac{2+\varepsilon}{3+\varepsilon}} (\int d^3x |{\bar \omega}({\bar x})|^{3+\varepsilon})^{\frac{1}{3+\varepsilon}} \\ \nonumber && \leq \sqrt{4 \pi} \sqrt{\int d^3 x |{\bar \nabla}{\bar \nabla} {\bar \omega}({\bar x})|^2} (\int d^3 x |{\bar \nabla}{\bar \omega}({\bar x})|^{3+\varepsilon})^{\frac{2+\varepsilon}{3+\varepsilon}} (\int d^3x |{\bar \omega} ({\bar x})|^{6})^{\frac{1+\varepsilon}{4(3+\varepsilon)}}(\int d^3x |\omega|^2({\bar x}) )^{\frac{3-\varepsilon}{4(3+\varepsilon)}} \\ \nonumber && \leq \sqrt{4 \pi} \sqrt{\int d^3 x |{\bar \nabla} {\bar \nabla} {\bar \omega}({\bar x})|^2} (\int d^3 x |{\bar \nabla} {\bar \omega}({\bar x})|^{3+\varepsilon})^{\frac{2+\varepsilon}{3+\varepsilon}} (\int d^3 x |{\bar \nabla}{\bar \omega} ({\bar x} )|^{2})^{\frac{1+\varepsilon}{4(1+\frac{\varepsilon}{3})}}
(\int d^3x |{\bar \omega}({\bar x})|^2)^{\frac{3-\varepsilon}{4(3+\varepsilon)}} \\ \nonumber && \leq \sqrt{4 \pi} \sqrt{\int d^3 x |{\bar \nabla}{\bar \nabla} {\bar \omega}({\bar x})|^2} (\int d^3 x |{\bar \nabla}{\bar \omega}({\bar x})|^{3+\varepsilon})^{\frac{2+\varepsilon}{3+\varepsilon}} (\int d^3 x |{\bar \nabla}{\bar \omega}({\bar x})|^2)^{\frac{1+\varepsilon}{4(1+\frac{\varepsilon}{3})}}(\int d^3 x |{\bar \omega}({\bar x})|^2)^{\frac{3-\varepsilon}{4(3+\varepsilon)}} \\ \nonumber && \leq \sqrt{4 \pi} (\int d^3 x |{\bar \nabla}{\bar \nabla}{\bar \omega}|^2({\bar x}) + 1)(\int d^3 x |{\bar \nabla}{\bar \omega}({\bar x})|^{3+\varepsilon}+1) (\int d^3 x |{\bar \nabla}{\bar \omega}({\bar x})|^2  +1 ) (\int d^3 x |{\bar \omega}({\bar x})|^2 +1).
\end{eqnarray}
The fourth integral in (\ref{ns30}) can be estimated as 
\begin{eqnarray}
&&
\int \int_{|y| \leq 1} d^3 x d^3 y |{\bar \nabla}{\bar \omega} ({\bar x})|^{2+\varepsilon} 
|{\bar \omega}({\bar x}+{\bar y})| \frac{|{\bar \nabla}{\bar \nabla}{\bar \omega}({\bar x})|}{|y|^2} \\ \nonumber &&= \int d^3 x |{\bar \nabla}{\bar \nabla}{\bar \omega}({\bar x})| |{\bar \nabla}{\bar \omega}({\bar x})|^{2+\varepsilon} (\int_{|y| \leq 1} d^3 y \frac{|{\bar \omega}({\bar x}+{\bar y})|}{|y|^2}) \\ \nonumber && \leq \int d^3 x |{\bar \nabla}{\bar \nabla}{\bar \omega}({\bar x})| |{\bar \nabla}{\bar \omega} ({\bar x})|^{2+\varepsilon} (\int d^3 y |{\bar \omega}({\bar x}+{\bar y})|^{6})^{\frac{1}{6}} (\int_{|y| \leq 1} \frac{d^3 y}{|y|^{\frac{12}{5}}})^{\frac{5}{6}} \\ \nonumber && \leq \frac{\nu}{25} \int d^3 x |{\bar \nabla}{\bar \nabla}{\bar \omega}({\bar x})|^2 |{\bar \nabla} {\bar \omega}({\bar x})|^{1+\varepsilon}+\frac{25}{4 \nu} (\frac{20 \pi}{3})^{\frac{5}{3}} \int d^3x |{\bar \nabla}{\bar \omega}({\bar x})|^{3+\varepsilon} (\int d^3 y |{\bar \nabla}{\bar \omega} ({\bar y})|^{2}) \\ \nonumber && \leq \frac{\nu}{25} \int d^3x |{\bar \nabla}{\bar \nabla}{\bar \omega}({\bar x})|^2 |{\bar \nabla}{\bar \omega}({\bar x})|^{1+\varepsilon} + \frac{25}{4 \nu}(\frac{20\pi}{3})^{\frac{5}{3}} \int d^3x |{\bar \nabla}{\bar \omega}({\bar x})|^{3+\varepsilon}  \int d^3 y |{\bar \nabla}{\bar \omega}({\bar y})|^2.
\end{eqnarray}
The fifth integral in (\ref{ns30}) can be estimated as
\begin{eqnarray}
&& \int d^3x \int_{|y| > 1} d^3y |{\bar \nabla}{\bar \omega}({\bar x})|^{2+\varepsilon} |{\bar \omega}({\bar x}+{\bar y})| \frac{|{\bar \nabla}{\bar \nabla}{\bar \omega}({\bar x})|}{|y|^2} \leq \\ \nonumber && \leq \int d^3x |{\bar \nabla}{\bar \omega}({\bar x})|^{2+\varepsilon} |{\bar \nabla}{\bar \nabla}{\bar \omega}({\bar x})| \int d^3y \frac{ |{\bar \omega}({\bar x}+{\bar y})|}{|y|^2} \leq \\ \nonumber && 
\leq \frac{\nu}{25} \int d^3x |{\bar \nabla}{\bar \nabla}{\bar \omega}({\bar x})|^2 |{\bar \nabla}{\bar \omega}({\bar x})|^{1+\varepsilon}  + \frac{25 \pi}{ \nu} \int d^3x |{\bar \nabla}{\bar \omega}({\bar x})|^{3+\varepsilon} (\int d^3y |{\bar \omega}({\bar x}+{\bar y})|^2 ) \\ \nonumber && \leq \frac{\nu}{25} \int d^3x |{\bar \nabla}{\bar \nabla}{\bar \omega}({\bar x})|^2 |{\bar \nabla}{\bar \omega}({\bar x})|^{1+\varepsilon} + \frac{25 \pi}{\nu} \int d^3x |{\bar \nabla}{\bar \omega}({\bar x})|^{3+\varepsilon} (\int d^3y |{\bar \omega}({\bar y})|^2).
\end{eqnarray}
Putting all the estimates together we obtain
\begin{eqnarray}\label{3+est}
&& \frac{1}{3+\varepsilon} \partial_t \int |{\bar \nabla}{\bar \omega}({\bar x})|^{3+\varepsilon} d^3x \leq \\ \nonumber && 
\leq 240 \pi \int |{\bar \nabla}{\bar \omega}({\bar x})|^{3+\varepsilon} d^3x  (\int |{\bar \nabla}{\bar \nabla} {\bar \omega}({\bar x})|^2 d^3x + 1) \\ \nonumber && \quad \quad \quad + 12 \sqrt{12 \pi} \int d^3 x |{\bar \nabla}{\bar \omega}({\bar x})|^{3+\varepsilon} (\int d^3y |{\bar \nabla}{\bar \omega}({\bar x}+{\bar y})|^2+1) + \\ \nonumber && + 6.6 \cdot 10^3 \sqrt{3} (\int |{\bar \nabla} {\bar \nabla} {\bar \omega}({\bar x})|^2 d^3x + 1)\int |{\bar \nabla}{\bar \omega}({\bar x})|^{3+\varepsilon} d^3x( \int |{\bar \nabla}{\bar \omega}({\bar x})|^2 d^3x +1)(N+1)+ \\ \nonumber && \quad \quad \quad  +
6 \sqrt{12 \pi} (\int |{\bar \nabla}{\bar \nabla} {\bar \omega}({\bar x})|^2 d^3x +1)\int |{\bar \nabla}{\bar \omega}({\bar x})|^{3+\varepsilon} d^3x (\int |{\bar \nabla}{\bar \omega}({\bar x})|^2 d^3x+1)(N+1) +\\ \nonumber && \quad \quad \quad \quad \quad + \frac{150 \sqrt{3}}{4 \nu}(\frac{20 \pi}{3})^{\frac{5}{3}} \int |{\bar \nabla}{\bar \omega}({\bar x})|^{3+\varepsilon} d^3x \int |{\bar \nabla}{\bar \omega}({\bar x})|^2 d^3x + \\ \nonumber && \quad \quad \quad \quad \quad \quad \quad  + \frac{150 \sqrt{3} \pi}{ \nu} \int |{\bar \nabla}{\bar \omega}({\bar x})|^{3+\varepsilon} d^3x (\int |{\bar \omega}({\bar y})|^2 d^3y) \\ \nonumber && + (6.6 \cdot 10^3 \sqrt{3} + 6 \sqrt{12 \pi})(\int d^3 x |{\bar \nabla}{\bar \nabla}{\bar \omega}({\bar x})|^2 + 1)(\int d^3 x |{\bar \nabla}{\bar \omega}({\bar x})|^2 +1) (N+1).
\end{eqnarray}
Estimating the right hand side of (\ref{3+est}) from above and integrating over $t$ we obtain 
\begin{eqnarray}\label{omega3.01estimate}
&& \int |{\bar \nabla}{\bar \omega}({\bar x})|^{3+\varepsilon} d^3x (t_2) \leq   \\ \nonumber &&\leq \left(\int |{\bar \nabla}{\bar \omega}({\bar x})|^{3+\varepsilon} d^3x (t_1) +\right. \\ \nonumber && \left.+ 3.55 \cdot 10^4 (N+1)( {\rm max}_{t \in [t_1, t_2]}( \int |{\bar \nabla}{\bar \omega}({\bar x})|^2 d^3x (t)+1) (\int_{t_1}^{t_2} dt \int |{\bar \nabla}{\bar \nabla} {\bar \omega}({\bar x})|^2 d^3x (t)+t_2-t_1)\right) \\  \nonumber &&
\cdot {\rm exp}\left(\frac{2.48 \cdot 10^3}{\nu} C(t_1, t_2) +  \right. \\ \nonumber && \left.\  \ \ + (2.28 \cdot 10^3 + 3.55 \cdot 10^4(N+1) ({\rm max}_{t \in [t_1, t_2]} \int |{\bar \nabla}{\bar \omega}({\bar x})|^2 d^3x+1)) \int_{t_1}^{t_2} dt \int |{\bar \nabla}{\bar \nabla}{\bar \omega}({\bar x})|^2 d^3 x\right. + \\  \nonumber && \quad \quad \quad \quad \quad \left. + (3.58 \cdot 10^4 + 3.55 \cdot 10^4 N +\frac{3.376 \cdot 10^4}{\nu}){\rm max}_{t \in [t_1, t_2]}\int |{\bar \nabla}{\bar \omega}({\bar x})|^{2} d^3x (t_2-t_1) +\right. \\ \nonumber && \quad \quad \quad \quad \quad \quad \quad \left.+(3.58 \cdot 10^4+3.55 \cdot 10^4 N)(t_2 -t_1) \right)\\ \nonumber && \leq (\int |{\bar \nabla}{\bar \omega}({\bar x})|^{3+\varepsilon}d^3x (t_1) +1) \\ \nonumber && {\rm exp}\left(\frac{2.48 \cdot 10^3}{\nu}C(t_1, t_2)+3.58 \cdot 10^4 (N+1)(t_2-t_1)+\right. \\ \nonumber && \left.\quad \quad \quad + 3.78 \cdot 10^4 (1+\int|{\bar \nabla}{\bar \omega}({\bar x})|^2 d^3x(t_1))^2 (N+1+\frac{1}{\nu})(t_2-t_1+ \frac{1.945}{\nu}) \right. \\ \nonumber && \left.\quad \quad \quad  \times {\rm exp}\left(2.45 \cdot 10^5\frac{(t_2-t_1)}{\nu}  + \frac{180}{\nu} C(t_1,t_2) + \frac{7.294 \cdot 10^7}{\nu} e^{4.22 N^{\frac{1}{2}} (t_2 - t_1)}\times  \right. \right. \\ \nonumber && \left. \left. \left(1+ \frac{3.99 \cdot 10^{26}}{\nu^{4.02}} N^{2.01} C(t_1, t_2) + \frac{3.562 \cdot 10^5}{\nu} N^{\frac{1}{2}} C(t_1, t_2) + 1.7 \cdot 10^5\frac{N}{\nu}+0.371 \cdot 10^5 (t_2-t_1)\right. \right. \right. \\ \nonumber && \qquad \qquad \qquad \qquad  \left. \left. \left. + 1.066 N^{\frac{1}{2}}(t_2 - t_1)+ |\omega|_{2+\frac{2}{3}\varepsilon}^{(2+\frac{2}{3}\varepsilon)(1-\alpha)}(t_1)\right)^5 \right) \right)
\end{eqnarray}

{\bf Part 5:} Taking the scalar product of $|\omega({\bar x})|^{n-2} {\bar \omega}({\bar x})$ with the Navier-Stokes equation for $\omega$ and integrating over the whole space we obtain
\begin{eqnarray}\label{eqn38}
&&\frac{1}{n} \partial_t \int |\omega({\bar x})|^n  d^3 x - \nu \int |\omega({\bar x})|^{n-2} {\bar \omega}({\bar x}) \cdot \Delta {\bar \omega}({\bar x}) d^3 x = \\ \nonumber && \qquad \qquad \qquad = \frac{3}{4\pi}\int (\int_{|y| \leq |{\omega}|_2^{-1}}+\int_{|y| > |{\omega}|_2^{-1}}) (|\omega|^{n-2}({\bar x}) {\bar \omega}({\bar x}) \cdot {\hat y}) ({\bar \omega}({\bar x}) \times {\bar \omega}({\bar x} + {\bar y}) \cdot {\hat y}) \frac{d^3 x d^3 y}{|y|^3}.
\end{eqnarray}
We rewrite the absolute value of the integral over $|y| \leq |{\bar \omega}|_2^{-1}$ on the right hand side of (\ref{eqn38}) as
\begin{eqnarray}
&&\left|\int d^3 x \int_{|y| \leq |{\omega}|_2^{-1}} d^3 y (|\omega|^{n-2} ({\bar x}) {\bar \omega}({\bar x}) \cdot {\hat y}) ({\bar \omega}({\bar x}) \times ({\bar \omega}({\bar x}+{\bar y}) - {\bar \omega}({\bar x})) \cdot {\hat y}) \frac{1}{|y|^3}\right| \qquad \qquad \qquad  \\ \nonumber && \leq 3 \sqrt{6} \int d^3 x \int_{|y| \leq |{\omega}|_2^{-1}} d^3 y \int_{0}^{1} d s |{\bar \omega}({\bar x})|^{n} \frac{|{\bar \nabla}{\bar \omega}({\bar x} + s {\bar y})|}{|y|^2} \qquad \qquad \qquad \\ \nonumber && \leq
3 \sqrt{6}\left(\frac{4 \pi(2+\varepsilon)}{\varepsilon}\right)^{\frac{2+\varepsilon}{3+\varepsilon}} \int d^3 x |{\bar \omega}({\bar x})|^n \int_{0}^{1} ds \left(\int_{|y| \leq |{\omega}|_2^{-1}} d^3 y  |{\bar \nabla}{\bar \omega}({\bar x}+ s{\bar y})|^{3+\varepsilon}\right)^{\frac{1}{3+\varepsilon}} \\ \nonumber && \leq 3 \sqrt{6}\left(\frac{4\pi(2+\varepsilon)}{\varepsilon}\right)^{\frac{2+\varepsilon}{3+\varepsilon}} \int d^3 x |{\bar \omega}({\bar x})|^n \int_{0}^{1} \frac{ds}{s^{\frac{3}{3+\varepsilon}}} (\int_{|z| \leq s|\omega|_2^{-1}} d^3 z |{\bar \nabla}{\bar \omega}({\bar x}+{\bar z})|^{3+\varepsilon})^{\frac{1}{3+\varepsilon}} \\ \nonumber && \leq 3 \sqrt{6} \frac{3+\varepsilon}{\varepsilon}\left(\frac{4 \pi(2 + \varepsilon)}{\varepsilon}\right)^{\frac{2+\varepsilon}{3+\varepsilon}}  \int d^3 x |{\bar \omega}({\bar x})|^{n} (\int |{\bar \nabla}{\bar \omega}({\bar x})|^{3+\varepsilon} d^3x)^{\frac{1}{3+\varepsilon}}\\ \nonumber && \leq \frac{3 \sqrt{6}(3+\varepsilon)}{\varepsilon}\left(\frac{4 \pi(2+\varepsilon)}{\varepsilon}\right)^{\frac{2+\varepsilon}{3+\varepsilon}} \int|{\bar \omega}({\bar x})|^n d^3x (\int |{\bar \nabla}{\bar \omega}({\bar x})|^{3+\varepsilon} d^3x)^{\frac{1}{3+\varepsilon}} \\ \nonumber && \leq 4.2 \cdot 10^5 \int |{\bar \omega}({\bar x})|^n d^3x \left(\int |{\bar \nabla}{\bar \omega}({\bar x})|^{3+\varepsilon} d^3 x \right)^{\frac{1}{3+\varepsilon}}. 
\end{eqnarray}
The absolute value  of the integral over $|y| > |{\bar \omega}|_2^{-1}$ can be estimated as
\begin{equation}
3 \sqrt{2}\int d^3 x |{\bar \omega}({\bar x})|^n \int_{|y| > |{\omega}|_2^{-1}} \frac{|{\bar \omega}({\bar x}+{\bar y})|}{|y|^3} \leq 3 \sqrt{8 \pi} N^2 \int d^3 x |{\bar \omega}({\bar x})|^n.
\end{equation}

Therefore we obtain from (\ref{eqn38}) that 
\begin{eqnarray}
&&\frac{1}{n} \partial_t \int |\omega|^n ({\bar x}) d^3 x \leq 1.1 \cdot 10^5 (N^2+1) \int |\omega|^n ({\bar x}) d^3 x\  {\rm max}_{t \in [t_1, t_2]} \left(\int |{\bar \nabla}{\bar \omega}({\bar x})|^{3+\varepsilon} d^3x (t)+1\right) \qquad 
\end{eqnarray}
We obtain that
\begin{eqnarray}
&& {\rm ln} \int |\omega|^n({\bar x}) d^3 x (t_2) \leq {\rm ln} \int |\omega|^n({\bar x}) d^3 x (t_1) + \\ \nonumber && \qquad \qquad \qquad + 1.1 \cdot 10^5 n (N^2+1) (t_2-t_1)\ {\rm max}_{t \in [t_1, t_2]} \left(\int |{\bar \nabla}{\bar \omega}({\bar x})|^{3+\varepsilon} d^3x (t)+ 1\right) 
\end{eqnarray}
Therefore
\begin{eqnarray}
&& \int |\omega|^n ({\bar x}) d^3 x (t_2) \leq \int |\omega|^n ({\bar x}) d^3 x (t_1)\\ \nonumber && \qquad \qquad \qquad  {\rm exp} \left(1.1 \cdot 10^5 n (N^2+1) (t_2-t_1){\rm max}_{t \in [t_1, t_2]} \left(|{\bar \nabla}{\bar \omega}|_{3+\varepsilon}^{3+\varepsilon}(t)+1\right) \right)\nonumber 
\end{eqnarray}
This implies
\begin{eqnarray}
&& |\omega|_n(t_2) \leq |\omega|_n(t_1)\ {\rm exp}\left(1.1 \cdot 10^5 (N^2+1) (t_2 -t_1) {\rm max}_{t \in [t_1, t_2]} \left(|{\bar \nabla}{\bar \omega}|_{3+\varepsilon}^{3+\varepsilon}(t) +1\right)\right) \nonumber 
\end{eqnarray}
Taking the limit $n \rightarrow \infty$ we obtain
\begin{eqnarray}\label{omegainftyupper}
&& |\omega|_{\infty}(t_2) \leq |\omega|_{\infty}(t_1)\ {\rm exp} \left(1.1 \cdot 10^5 (N^2+1) (t_2 - t_1) \ {\rm max}_{t \in [t_1, t_2]} \left(|{\bar \nabla}{\bar \omega}|_{3+\varepsilon}^{3+\varepsilon}(t)+1\right)\right)\quad \quad  \quad  \\ \nonumber 
&&  \leq |\omega|_{\infty}(t_1) {\rm exp}\left(1.1 \cdot 10^5 (N^2 + 1)(t_2 - t_1) (|{\bar \nabla}{\bar \omega}|_{3.01}^{3.01}(t_1)+2) \times \right. \\ \nonumber && \left. \times  {\rm exp}\left(\frac{2.48 \cdot 10^3}{\nu}C(t_1,t_2)+3.58\cdot 10^4(N+1)(t_2-t_1)+ \right. \right. \\ \nonumber && \left. \left.\quad \quad \quad  + 3.78 \cdot 10^4(1+ |{\bar \nabla}{\bar \omega}|_{2}^{2}(t_1))^2 (N+1+\frac{1}{\nu})(t_2-t_1+\frac{1.945}{\nu})  \right. \right. \\ \nonumber && \quad \quad \quad   \left. \left. \times  {\rm exp}\left( 2.45 \cdot 10^5 \frac{(t_2 - t_1)}{\nu} + \frac{180}{\nu} C(t_1, t_2) + \frac{7.294 \cdot 10^7}{\nu} e^{4.22 N^{\frac{1}{2}}(t_2-t_1)} \times  \right. \right. \right. \\ \nonumber && \qquad \qquad \qquad \left. \left. \left. \times \left(1+ \frac{3.99 \cdot 10^{26}}{\nu^{4.02}} N^{2.01} C(t_1, t_2) + \frac{3.562 \cdot 10^5}{\nu} N^{\frac{1}{2}} C(t_1, t_2) + 1.7 \cdot 10^5 \frac{N}{\nu}+ \right. \right. \right. \right. \\ \nonumber && \qquad \qquad \qquad \left. \left. \left. \left.+0.371 \cdot 10^5(t_2-t_1) + 1.066 N^{\frac{1}{2}}(t_2-t_1) + |{\bar \nabla}{\bar \omega}|_{3.01}^{0.76153}(t_1)\right)^5\right) \right) \right).
\end{eqnarray}  
As a result we obtain that if $|\omega|_2 (t) \leq N$ for $t \in [0, T- \delta_N]$, then for $\delta \in (0, T- \delta_N)$
\begin{eqnarray}\label{omegainftyestN}
&& |\omega|_{\infty} (t) \leq |\omega|_{\infty}(\delta) \ {\rm exp} \left(1.1 \cdot 10^5 (N^2+1) t \ {\rm max}_{t \in [\delta, T- \delta_N]} \left( |{\bar \nabla}{\bar \omega}|_{3+\varepsilon}^{3+\varepsilon} (t) + 1 \right) \right) \quad \quad \quad \\ \nonumber 
&& \leq |\omega|_{\infty}(\delta) {\rm exp}\left(1.1 \cdot 10^5 (N^2+1) t (|{\bar \nabla}{\bar \omega}|_{3.01}^{3.01}(\delta)+2) \times \right. \\ \nonumber && \left. \times {\rm exp} \left(\frac{2.48 \cdot 10^3}{\nu}C(t)+3.58 \cdot 10^4 (N+1)t+
\right. \right. \\ \nonumber && \left. \left.  \quad \quad \quad  +3.78 \cdot 10^4 (1+|{\bar \nabla}{\bar \omega}|_{2}^{2}(\delta))^2 (N+1+\frac{1}{\nu})(t+\frac{1.945}{\nu}) \right. \right. \\ \nonumber && \left. \left. \quad \quad \quad  {\rm exp} \left( 2.45 \cdot
10^5 \frac{t}{\nu} + \frac{180}{\nu} C(t) + \frac{7.294 \cdot 10^7}{\nu} e^{4.22 N^{\frac{1}{2}} t}  \right. \right. \right. \\ \nonumber && \left. \left. 
\left. \qquad \qquad \qquad \times \left(1+\frac{3.99 \cdot 10^{26}}{\nu^{4.02}} N^{2.01} C(t) + \frac{3.562 \cdot 10^5}{\nu}
N^{\frac{1}{2}} C(t) + 1.7 \cdot 10^5 \frac{N}{\nu}) + \right. \right. \right. \right. \\ \nonumber && \qquad \qquad \qquad \left. \left. \left. \left. +0.371 \cdot 10^5 t+1.066 N^{\frac{1}{2}} t+ |{\bar \nabla}{\bar \omega}|_{3.01}^{0.76153}(\delta) \right)^5 \right) \right) \right).
\end{eqnarray}
Now we shall improve the estimate on $|\omega|_{\infty}(t)$.

\begin{center}{\bf Proof of Theorem 2.} \end{center}

\noindent {\bf Part 1:} Suppose first that there exists $n_0$ such that for all $n > n_0$, $|\omega|_n(t) \geq e^e$ for all $t \in [t_1, t_2]$. Let $n > n_0$.

As in the derivation of (\ref{eqn38}) we obtain
\begin{eqnarray}\label{eq1omn}
&&\frac{1}{n} \partial_t \int |\omega|^n d^3 x  -  \nu \int {\bar \nabla}[|\omega|^{n-2}({\bar x}) {\bar \omega}({\bar x})] ({\bar \nabla} {\bar \omega}({\bar x}))d^3 x  \\ \nonumber && \qquad \qquad \qquad  =\frac{3}{4 \pi} \int \int [|\omega|^{n-2}({\bar x}) {\bar \omega}({\bar x}) \cdot {\hat y}] [{\bar \omega}({\bar x}) \times {\bar \omega}({\bar x} + {\bar y}) \cdot {\hat y}] \frac{d^3 x d^3 y}{|y|^3}.
\end{eqnarray}
Since we know from the Theorem 1 that $|\omega|_{\infty}(t) < \infty$ for all $t \in [t_1, t_2]$, we can estimate the right hand side of (\ref{eq1omn}) as follows
\begin{eqnarray}
&& \frac{1}{n} \partial_t \int |\omega|^n d^3 x  \leq \frac{9 \sqrt{6}}{4 \pi}\int_{0}^{1} d s \int |\omega|^n ({\bar x}) d^3 x \int_{|y| \leq {\rm min}(|\omega|_{\infty}^{-1}, 1)} \frac{|{\bar \nabla} {\bar \omega}|({\bar x} + s {\bar y}) }{|y|^2} d^3 y\nonumber \\ &&  + \frac{9 \sqrt{2}}{4 \pi}\int |\omega|^n ({\bar x}) d^3 x \int_{1 > |y| \geq {\rm min}(|\omega|_{\infty}^{-1}, 1)} \frac{|\omega|({\bar x} + {\bar y})}{|y|^3 } d^3 y + \frac{9 \sqrt{2}}{4 \pi} \int |\omega|^n({\bar x}) d^3x \int_{|y| \geq 1} \frac{|\omega|({\bar x}+{\bar y})}{|y|^3} d^3 y .\nonumber
\end{eqnarray} 
Therefore
\begin{eqnarray}
&& \frac{1}{n} \partial_t \int |\omega|^n ({\bar x}) d^3 x \leq 9 \sqrt{6}\ \frac{3+\varepsilon}{4 \pi\varepsilon} \left(\frac{4 \pi(2+\varepsilon)}{\varepsilon}\right)^{\frac{2+\varepsilon}{3+\varepsilon}} \int |\omega|^n ({\bar x}) d^3 x  (\int |{\bar \nabla}{\bar \omega}({\bar x})|^{3+\varepsilon} d^3x)^{\frac{1}{3+\varepsilon}}+\nonumber \\ \nonumber && \quad  \quad  \quad   + 9 \sqrt{2}|\omega|_{\infty}(t) {\rm ln}({\rm max}( |\omega|_{\infty}(t), 1)) \int |\omega|^n ({\bar x}) d^3 x + 3 \sqrt{\frac{3 N}{2 \pi}} \int |\omega|^n ({\bar x}) d^3 x.
\end{eqnarray}

\begin{eqnarray}
&& \frac{1}{n} \partial_t \int |\omega|^n ({\bar x}) d^3 x \leq \int |\omega|^n ({\bar x}) d^3 x \cdot \nonumber  \\ \nonumber && \quad \quad \quad \cdot \left(12.78 |\omega|_{\infty}(t) {\rm ln} ({\rm max}(|\omega|_{\infty}(t), 1)) + 1.1 \cdot 10^5 |{\bar \nabla}{\bar \omega}|_{3+\varepsilon}(t) + 2.07 \sqrt{N}\right).\nonumber 
\end{eqnarray}

\begin{eqnarray}
  && \frac{\partial_t \int |\omega|^n ({\bar x}) d^3 x}{\int |\omega|^n ({\bar x}) d^3 x} \leq n \left( 12.78 |\omega|_{\infty}(t) {\rm ln}({\rm max} |\omega|_{\infty} (t), 1))  + 1.1 \cdot 10^5 |{\bar \nabla}{\bar \omega}|_{3+\varepsilon}(t) + 2.07 \sqrt{N}\right).\nonumber 
\end{eqnarray}

\begin{eqnarray}
&& \partial_t {\rm ln} \int |\omega|^n ({\bar x}) d^3 x \leq n \left(12.78 |\omega|_{\infty}(t) {\rm ln}({\rm max} (|\omega|_{\infty}(t), 1)) + 1.1 \cdot 10^5 |{\bar \nabla}{\bar \omega}|_{3+\varepsilon}(t) + 2.07 \sqrt{N} \right). \nonumber 
\end{eqnarray}

\begin{eqnarray}
&& \frac{\partial_t {\rm ln} \int |\omega|^n ({\bar x}) d^3 x}{{\rm ln} \int |\omega|^n ({\bar x}) d^3 x} \leq \frac{ n \left(12.78 |\omega|_{\infty}(t) {\rm ln} ({\rm max}(|\omega|_{\infty}(t),1))  + 1.1 \cdot 10^5 |{\bar \nabla}{\bar \omega}|_{3+\varepsilon}(t) + 2.07 \sqrt{N}\right)}{{\rm ln} \int |\omega|^n ({\bar x}) d^3 x}\nonumber
\end{eqnarray} 

\begin{eqnarray}
&& \frac{\partial_t {\rm ln~ln} \int |\omega|^n({\bar x}) d^3 x }{{\rm ln~ln} \int |\omega|^n ({\bar x}) d^3 x} \leq \frac{n \left(12.78 |\omega|_{\infty}(t) {\rm ln}({\rm max}( |\omega|_{\infty}(t),1))+ 1.1 \cdot 10^5 |{\bar \nabla}{\bar \omega}|_{3 + \varepsilon}(t) + 2.07 \sqrt{N}\right)}{{\rm ln} \int |\omega|^n ({\bar x}) d^3 x~ {\rm ln ~ln} \int |\omega|^n({\bar x}) d^3 x} \nonumber 
\end{eqnarray}                                                                                                                                                 
\begin{eqnarray}
&&\partial_t {\rm ln~ln~ln} \int |\omega|^n({\bar x}) d^3 x \leq \frac{n \left(12.78 |\omega|_{\infty}(t) {\rm ln} ({\rm max}(|\omega|_{\infty}(t),1))+ 1.1 \cdot 10^5 |{\bar \nabla}{\bar \omega}|_{3+\varepsilon}(t) + 2.07 \sqrt{N}\right)}{n {\rm ln} ((\int |\omega|^n ({\bar x}) d^3 x)^{\frac{1}{n}}) {\rm ln} (n {\rm ln}(\int |\omega|^n ({\bar x}) d^3 x)^{\frac{1}{n}})} \nonumber 
\end{eqnarray} 

\begin{eqnarray}
&& \partial_t {\rm ln~ln~ln} \int |\omega|^n({\bar x}) d^3 x \leq \frac{12.78 |\omega|_{\infty}(t) {\rm ln}({\rm max}( |\omega|_{\infty}(t),1))+ 1.1 \cdot 10^5 |{\bar \nabla}{\bar \omega}|_{3+\varepsilon}(t) + 2.07 \sqrt{N} }{{\rm ln} |\omega|_n(t) ({\rm ln} n + {\rm ln} {\rm ln} |\omega|_n (t))}\nonumber 
\end{eqnarray} 

Below we consider $n$ satisfying the following condition
\begin{eqnarray}\label{nlowerbound}
{\rm ln} n >{\rm max}_{t \in [t_1, t_2]} \left(12.78 |\omega|_{\infty}(t) {\rm ln}({\rm max}( |\omega|_{\infty}(t),1))+ 1.1 \cdot 10^5 |{\bar \nabla}{\bar \omega}|_{3+\varepsilon} (t) + 2.07 \sqrt{N}\right), \qquad 
\end{eqnarray}
in fact, we shall consider even larger $n$, because instead of $|\omega|_{\infty}(t)$ we shall use the estimate (\ref{omegainftyestN})
\begin{eqnarray}
|\omega|_{\infty}(t) \leq |\omega|_{\infty}(t_1) {\rm exp}\left(1.1 \cdot 10^5 (N^2 +1) (t_2-t_1) {\rm max}_{t \in [t_1, t_2]} \left(|{\bar \nabla}{\bar \omega}|_{3+\varepsilon}^{3+\varepsilon}(t)+1 \right)\right), \ \ \ \nonumber 
\end{eqnarray}
and instead of $|{\bar \nabla}{\bar \omega}|_{3+\varepsilon}^{3+\varepsilon}$ we shall use (\ref{omega3.01estimate})
\begin{eqnarray}
&& \int|{\bar \nabla}{\bar \omega}({\bar x})|^{3+\varepsilon} d^3x (t) \leq (\int |{\bar \nabla}{\bar \omega}({\bar x})|^{3+\varepsilon} d^3 x (t_1)+1) \times  \nonumber \\ \nonumber && \ \ \ \ \ \times {\rm exp} \left(\frac{2.48 \cdot 10^3}{\nu} C(t_1,t_2)+3.58 \cdot 10^4 (N+1) (t_2-t_1) + \right. \nonumber \\ \nonumber && \left. \quad \quad \quad  + 3.78 \cdot 10^4 (1+\int|{\bar \nabla}{\bar \omega}({\bar x})|^2 d^3 x(t_1))^2 (N+1+\frac{1}{\nu})(t_2-t_1+\frac{1.945}{\nu}) \right.  \\ \nonumber && \left. \quad \quad \quad \times {\rm exp}\left( \frac{2.45 \cdot 10^5}{\nu}  (t_2-t_1)  + \frac{180}{\nu} C(t_1, t_2) + \frac{7.294 \cdot 10^7}{\nu^2} e^{4.22 N^{\frac{1}{2}} (t_2-t_1)} \times \right. \right. \\ \nonumber && \left. \left. \qquad \qquad \qquad \left(1+\frac{3.99 \cdot 10^{26}}{\nu^{4.02}} N^{2.01} C(t_1, t_2) + \frac{3.562 \cdot 10^5}{\nu} N^{\frac{1}{2}} C(t_1, t_2) + 1.7 \cdot 10^5 \frac{N}{\nu}+ \right. \right. \right.  \\ \nonumber && \qquad \qquad \qquad \left. \left. \left. + 0.371 \cdot 10^5(t_2-t_1) + 1.066 N^{\frac{1}{2}} (t_2-t_1) + |\omega|_{\frac{6.02}{3}}^{0.76153}(t_1) \right)^5 \right) \right).
\end{eqnarray}
Therefore consider $n$ such that 
\begin{eqnarray}\label{nestimate}
&&  {\rm ln} n \geq 2.07 \sqrt{N} + 12.78 (|\omega|_{\infty}(t_1)+1) (|{\rm ln}(|\omega|_{\infty}(t_1))|+1)  \times \\  \nonumber && \ \ \ \ \  \times {\rm exp} \left(2.2 \cdot 10^5 (N^2 +1) (t_2-t_1+1) (|{\bar \nabla}{\bar \omega}({\bar x})|^{3+\varepsilon}_{3+\varepsilon}(t_1) + 1)  \right. \\ \nonumber && \quad \quad \quad \times \left. {\rm exp}\left( \frac{2.48 \cdot 10^3}{\nu} C(t_1, t_2)+ 3.58 \cdot 10^4 (N +1)(t_2-t_1)+ \right. \right. \\ \nonumber && \quad \quad \quad \left. \left. +3.78 \cdot 10^4(1+|{\bar \nabla}{\bar \omega}|_{2}^2(t_1))^2 (N+1+ \frac{1}{\nu})(t_2-t_1+\frac{1.945}{\nu}) \right. \right.  \\ \nonumber && \quad \quad \quad \left. \left. {\rm exp} \left( 2.45 \cdot 10^5 \frac{(t_2-t_1)}{\nu} + \frac{180}{\nu}C(t_1,t_2)+\frac{7.294 \cdot 10^7}{\nu^2} e^{4.22 N^{\frac{1}{2}}(t_2-t_1)} \right. \right. \right. \\ \nonumber && \qquad \qquad \qquad  +\left. \left.\left.  \left(1+\frac{3.99 \cdot 10^{26}}{\nu^{4.02}} N^{2.01} C(t_1, t_2)+ \frac{3.562 \cdot 10^5}{\nu} N^{\frac{1}{2}} C(t_1, t_2)+ 1.7 \cdot 10^5 \frac{N}{\nu}+ \right. \right. \right. \right. \\ \nonumber && \qquad \qquad  \qquad \left. \left. \left. \left.  + 0.371 \cdot 10^5 (t_2-t_1)+ 1.066 N^{\frac{1}{2}} (t_2-t_1) + |{\bar \nabla}{\bar \omega}|_{3+\varepsilon}^{0.76153}(t_1) \right)^5 \right) \right)\right)
\end{eqnarray}
   
For $n$ satisfying (\ref{nestimate}) we have that               
\begin{eqnarray}
\partial_t {\rm ln~ln~ln} \int |\omega|^n ({\bar x}) d^3 x \leq 1.
\end{eqnarray}
After integrating we obtain
\begin{eqnarray}
{\rm ln~ln~ln} \int |\omega|^n ({\bar x}, t) d^3 x \leq t-t_1 + {\rm ln~ln~ln} \int |\omega|^n ({\bar x}, t_1) d^3 x.
\end{eqnarray}
Therefore we obtain
\begin{eqnarray}
\int |\omega|^n ({\bar x}, t) d^3 x \leq e^{e^{e^{t-t_1}}} \int |\omega|^n ({\bar x}, t_1) d^3 x ~{\rm on} ~t\in [t_1, t_2].
\end{eqnarray}

Therefore we  obtain that for all $n$ satisfying (\ref{nestimate})
\begin{eqnarray}\label{omginftestimate}
|\omega|_n (t) \leq e^{\frac{1}{n} e^{e^(t-t_1)}} |\omega|_n (t_1) \leq e^{\frac{1}{n} e^{e^(t-t_1)}} |\omega|_{\infty}(t_1) ~{\rm ~for~} t \in [t_1, t_2].
\end{eqnarray}

Taking $n \rightarrow \infty$ we obtain that 
\begin{eqnarray}\label{omegainftestimate}
|\omega|_{\infty}(t) \leq |\omega|_{\infty}(t_1) ~{\rm for}~ t \in [t_1, t_2].
\end{eqnarray}

\noindent{\bf Part 2:} Suppose now that there exists $n_0 > 0$ such that for all $n > n_0$, $|\omega|_n(t)$ is a continuous function of $t$ for all $t \in [t_1, t_2]$. Fix $n > n_0$  such that the condition (\ref{nestimate}) is satisfied.

Suppose $|\omega|_n(t_1) \geq e^e$. By continuity of $|\omega|_n(t)$ we obtain that there exists a set consisting of a union of open intervals $\cup_j I_j$ such that on each interval $I_j$, $|\omega|_n(t) < e^e$ for all $t \in I_j$. Let $\left(\cup_j I_j\right)^{\rm c} = \cup_k I_k'$, where each interval $I_k' = [t_{1_k}', t_{2_k}']$ is such that $|\omega|_n(t) \geq e^e$ for all $t \in I_k'$ and $|\omega|_n(t_{1_k}') = e^e$ for $k > 1$ and $|\omega|_n(t_{1_k}') = |\omega|_n(t_1)$ for $k = 1$.

Then by formula (\ref{omginftestimate})  we obtain on each interval $I_k'$
\begin{eqnarray}
|\omega|_n(t) \leq e^{\frac{1}{n}e^{e^{t_{2_k}' - t_{1_k}'}}} e^e  ~{\rm for}~ k>1.
\end{eqnarray}

For the first interval $I_1' = [t_1, t_{2_1}']$ we obtain that 
\begin{eqnarray}
|\omega|_n(t) \leq e^{\frac{1}{n}e^{e^{t_{2_1}'- t_1}}} |\omega|_n(t_1).
\end{eqnarray}
Therefore 
\begin{eqnarray}
|\omega|_n(t) \leq e^{\frac{1}{n} e^{e^{t_2-t_1}}} |\omega|_n(t_1) ~{\rm for ~ all}~ t \in [t_1, t_2].
\end{eqnarray}

Suppose $|\omega|_n(t_1) < e^e.$ By continuity of $|\omega|_n(t)$ we obtain that there exists a set consisting of a union of open intervals $\cup_j I_j$ such that on each interval $I_j$, $|\omega|_n(t) < e^e$ for all $t \in I_j$.
Let $\left(\cup_j I_j\right)^{\rm c} = \cup_k I_k'$ where each interval $I_k' = [t_{1_k}', t_{2_k}']$ is such that $|\omega|_n(t) \geq e^e$ for all $t \in [t_{1_k}', t_{2_k}']$ and $|\omega|_n(t_{1_k}') = e^e$. 

Then by formula (\ref{omginftestimate}) we obtain that on each interval $I_k'$
\begin{eqnarray}
|\omega|_n(t) \leq e^{\frac{1}{n} e^{e^{t_{2_k}'-t_{1_k}'}}} |\omega|_n (t_{1_k}') ~{\rm for ~ all} ~ t \in [t_{1_k}', t_{2_k}'].
\end{eqnarray} 

Therefore if $|\omega|_n(t_1) < e^e$ we obtain
\begin{eqnarray}
|\omega|_n(t) \leq e^{\frac{1}{n} e^{e^{t_2-t_1}}} e^e ~ {\rm for~ all} ~ t \in [t_1, t_2].\end{eqnarray}

Therefore for every $n$ satisfying the condition (\ref{nestimate}) we obtain that
\begin{eqnarray}
|\omega|_n(t) \leq e^{\frac{1}{n}e^{e^{t-t_1}}} {\rm  max}(|\omega|_n(t_1), e^e) \leq e^{\frac{1}{n}e^{e^{t-t_1}}} {\rm max}(|\omega|_{\infty}(t_1), e^e), \quad {\rm for} \ t \in [t_1, t_2].\ \ \ 
\end{eqnarray} 
Taking the limit $n \rightarrow \infty$ we obtain that
\begin{eqnarray}
|\omega|_{\infty}(t) \leq {\rm max}(|\omega|_{\infty}(t_1), e^e) ~{\rm for~ all}~ t \in [t_1, t_2].
\end{eqnarray}

\begin{center}{\bf Proof of Theorem 3.}\end{center}

{\it {\bf Theorem \ 3:} Suppose $|\omega|_{\infty}(0) < \infty$. Then for any $T > 0$ there exists a unique solution of the 3D Navier Stokes Equation on $[0, T]$ such that $|\omega|_{\infty}(t) \leq {\rm max}(|\omega|_{\infty}(\delta), e^e)< \infty$ for $t \in [0, T]$ with $\delta = \frac{\nu}{2 C |\omega|_{\infty}^2(0)}$}, (where $C$ is a constant independent of $\omega$, $\nu$, etc.)

{\it {\bf Proof:}}  
First we shall construct an argument to show that if certain assumptions on $\omega$ hold then the solution is bounded. This argument will be used below.

Suppose we have that $|\omega|_{\infty}(0) <  \infty$. Then by the Local Existence and Uniqueness Theorem (see Theorem 4.2 in [Ku], where we take $p=2$ and $M_{2p} = |\omega|_{\infty}(0)$) there exists $t_0 = \frac{\nu}{2 C |\omega|_{\infty}^2(0)}$ such that for $t \in [0, t_0]$, we have that there exists a unique solution such that $|\omega|_2(t) < \infty$ and also by the analyticity results in the same theorem we know that for all $\delta \in (0, t_0]$, $|\omega|_{\infty}(\delta) < \infty$, $|\omega|_{3+\varepsilon}(\delta) < \infty$, $|{\bar \nabla}{\bar \omega}|_2 (\delta) < \infty$, $|{\bar \nabla}{\bar \omega}|_{3+\varepsilon}(\delta) < \infty$. Let $\delta = t_0.$ Assume that $|\omega|_2(t) \rightarrow \infty$ as $t \rightarrow T$ and $|\omega|_2(t) < \infty$ for $t \in [0, T)$. Also assume that for all $n$ large enough ($n$ satisfying (\ref{nestimate}) with $t_1 = \delta$), $|\omega|_n(t)$ is continuous as a function of $t$ for all $t \in [0, T)$. Then for any $N > N_0$ (for sufficiently large $N_0$) there exists $\delta_N > 0$ such that $|\omega|_2(t) \leq N$, for all $t \in [0, T-\delta_N]$. Then by Theorem 1 and Theorem 2 we obtain that $|\omega|_2(t) \leq |\omega|_{\infty}(t) \leq {\rm max}(|\omega|_{\infty}(\delta), e^e)$ for $t \in [\delta,T-\delta_N]$. By taking $N \rightarrow \infty$ we obtain that $|\omega|_2(t) \leq |\omega|_{\infty}(t) \leq {\rm max}(|\omega|_{\infty}(\delta), e^e)$ for $t \in [\delta, T)$. 
By Local Existence and Uniqueness Theorem (see Theorem 4.2 in [Ku]), we obtain that there exists $t_0 > 0$ such that there exists a unique solution on $t \in [T-\frac{t_0}{2}, T+\frac{t_0}{2}]$ that is continuous and analytic and such that $|\omega|_{\infty}(t) < \infty$ for $t \in [T-\frac{t_0}{2}, T+\frac{t_0}{2}]$ and then by Theorem 1 and Theorem 2 we obtain that 
\begin{eqnarray}
  |\omega|_2 (t) \leq |\omega|_{\infty}(t) \leq {\rm max}(|{\omega}|_{\infty}(\delta), e^e) \ {\rm for}\ t \in [\delta, T + \frac{t_0}{2}].
\end{eqnarray} 

Thus we showed that if certain assumptions on $\omega$ hold  then for any $T > 0$ the solution ${\bar \omega}({\bar x}, t)$ is finite for all ${\bar x} \in {\bf R}^3$ and $t \in [0, T]$.

To show existence and uniqueness of a solution on $[0, T]$, we shall repeatedly use the Local Existence Uniqueness Theorem (Theorem 4.2 in [Ku]).

Since $|\omega|_{\infty}(0) < \infty$ there exits $t_0 > 0$ given by $ t_0 = \frac{\nu}{2 C |\omega|_{\infty}^2(0)}$ such that there exists a unique solution of 3D Navier Stokes Equations that is continuous and analytic on $[0, t_0]$.  Let $\delta = t_0$, then $|\omega|_{\infty}(\delta) < \infty$, $|{\bar \nabla}{\bar \omega}|_{3+\varepsilon}(\delta) < \infty$. 

Now let $t_0' = \frac{\nu}{2 C {\rm max}(|\omega|_{\infty}^2(\delta), e^{2e})}$ and consider $t \in [\delta, \delta+ t_0']$. Using ${\bar \omega}({\bar x}, \delta)$ as the initial condition, we obtain that by the Local Existence and Uniqueness results there exists a unique solution which is continuous and analytic for $t \in [\delta, \delta+t_0']$. By Theorem 1, Theorem 2 and arguments above we obtain that  $|\omega|_2(t) \leq |\omega|_{\infty}(t) \leq {\rm max}(|\omega|_{\infty}(\delta), e^e)$ for all $t \in [\delta, \delta+ t_0']$.

Repeating this argument $k$ times, we obtain that for $t \in [\delta+(k-1)t_0', \delta+ k t_0']$ there exists a unique solution of the 3D Navier Stokes Equation such that $|\omega|_2(t) \leq {\rm max}(|\omega|_{\infty}(\delta+(k-1)t_0'), e^e) \leq {\rm max}(|\omega|_{\infty}(\delta), e^e)$ and $|\omega|_{\infty}(t) \leq {\rm max}(|\omega|_{\infty}(\delta+(k-1)t_0'), e^e) \leq {\rm max}(|\omega|_{\infty}(\delta), e^e)$.

For any $T>0$, we can continue this argument up to $k$ such that $k \geq \frac{T}{t_0'}$, to obtain that there exists a unique solution of the 3D Navier Stokes Equation on $[0,T]$ such that,
\begin{eqnarray}
|\omega|_2(t) \leq {\rm max}(|\omega|_{\infty}(\delta), e^e) \ {\rm and} \  |\omega|_{\infty}(t) \leq {\rm max}(|\omega|_{\infty}(\delta), e^e) \ {\rm for} \ t \in [\delta,T].
\end{eqnarray} 

\begin{center}{\bf Proof of Theorem 4.}\end{center}

{\it {\bf Theorem \ 4:} Suppose that $|\omega|_4(0) < \infty$. Then for any $T > 0$ there exists a unique solution of the 3D Navier Stokes Equation on $[0, T]$.}

{\it {\bf Proof:}} 
First we construct the argument that shows that if certain assumptions on $\omega$ hold, then the solution remains bounded. This argument is used later in the proof.

Suppose we have that $|\omega|_{4}(0)  < \infty$. Then by the Local Existence and Uniqueness Theorem (Theorem 4.2 in [Ku]) there exists $t_0' = \frac{\nu}{2 C |\omega|_{4}^2(0)}$ such that for $t \in [0, t_0']$, there exists a unique solution which is continuous on $(0, t_0']$. Therefore $|\omega|_2(t) < \infty$ for $t \in [0, t_0']$. In addition, by analyticity results in Theorem 4.2 of [Ku] we obtain that for any $\delta \in (0, t_0']$, $|\omega|_{\infty}(\delta) < \infty$, and also we have that $|\omega|_{3.01} (\delta) < \infty$ and $|{\bar \nabla}{\bar \omega}|_{3.01} (\delta) < \infty$. Suppose that $\delta = t_0'$.

Suppose that for some $T > 0$, $|\omega|_2(t) \rightarrow_{t \rightarrow T} \infty$, but $|\omega|_2 (t) < \infty$ for $t \in [0, T)$. (By the argument above we know that $T > t_0'$). Suppose also that there exists $n_0$ (satisfying (\ref{nestimate}) with $t_1= \delta$) such that for all $n > n_0$, $|\omega|_n(t)$ is a continuous function of $t$ for all $t \in [0, T)$.
For any $N > N_0$ (for a sufficiently large $N_0$), there exists $\delta_N > 0$ such that $|\omega|_2(t) \leq N$, for all $t \in [\delta, T-\delta_N]$. Then by Theorem 1, Theorem 2 and taking the limit $N \rightarrow \infty$, we obtain that $|\omega|_2(t) \leq |\omega|_{\infty}(t) \leq {\rm max}(|\omega|_{\infty}(\delta), e^e)$ for all $t \in [\delta, T)$. By the Local Existence and Uniqueness Theorem (Theorem 4.2 in [Ku]) there exists $t_0'' > 0$ such that there exists a unique solution which is continuous and analytic on $[T-\frac{t_0''}{2}, T+\frac{t_0''}{2}]$. Thus for any $T>0$ we obtain that $|\omega|_2(t) \leq |\omega|_{\infty}(t) \leq {\rm max}(|\omega|_{\infty}(\delta), e^e)$ for $t \in [\delta, T]$.

To show existence and uniqueness on $[0, T]$, we shall repeatedly use the Local Existence Uniqueness Theorem. Consider ${\bar \omega}({\bar x}, \delta)$ as the initial condition. Then there exists $t_0'' = \frac{\nu}{2 C {\rm max}(|\omega|_{\infty}^2(\delta), e^{2e})}$ such that there exists a unique solution which is continuous and analytic  on $[\delta, t_0'']$. By Theorem 1, Theorem 2 and the arguments above we have that $|\omega|_2(t) \leq |\omega|_{\infty}(t) \leq {\rm max}(|\omega|_{\infty}(\delta), e^e)$ for $t \in [\delta, t_0'']$. Let $t \in [t_0'', 2 t_0'']$. Then by the Local Existence and Uniqueness Theorem (see Theorem 4.2 in [Ku]) there exists a unique solution which is continuous and analytic  on $[t_0'', 2 t_0'']$. By Theorem 1, Theorem 2 and the arguments above this solution satisfies $|\omega|_2(t) \leq |\omega|_{\infty}(t) \leq {\rm max}(|\omega|_{\infty}(t_0''), e^e) \leq {\rm max}(|\omega|_{\infty}(\delta), e^e)$. Repeating the argument $k$ times we obtain that there exists a unique solution which is continuous and analytic on $[(k-1)t_0'', k t_0'']$, such that $|\omega|_2 (t) \leq |\omega|_{\infty}(t) \leq {\rm max}(|\omega|_{\infty}((k-1)t_0''), e^e) \leq {\rm max}(|\omega|_{\infty}(\delta), e^e)$ for $t \in [(k-1) t_0'', k t_0'']$.
Thus for any $T > 0$, letting $k \geq \frac{T}{t_0''}$ we obtain that there exists a unique solution of 3D Navier Stokes Equations on $[0, T]$ such that $|\omega|_2(t) \leq |\omega|_{\infty}(t) \leq {\rm max}(|\omega|_{\infty}(\delta), e^e)$ for $t \in [0, T]$.  

{\bf Bibliography}
\begin{itemize}
\item{}[CF] P. Constantin, C. Fefferman, Direction of vorticity and the problem of global regularity for the Navier-Stokes equations, {\it IUMJ}, {\bf Vol. 42}, (1993), 775-789.
\item{}[FMRT] Foias, Manley, Rosa, Temam, Navier-Stokes Equations and Turbulence, {\it Encyclopedia of Mathematics}, Cambridge Univ Press, (2001).
\item{}[G] Grujic, The geometric structure of the super-level sets and regularity for 3D Navier-Stokes equations, {\it IUMJ}, {\bf Vol. 50}, (2001), 1309-1317.
\item{}[GR] Grujic, Ruzmaikina, Interpolation between algebraic and geometric conditions for smoothness of the vorticity of the 3-D NSE, {\it IUMJ}, {\bf Vol. 53}, {\bf No. 4}, (2004), 1073-1080.
\item{}[J] J. D. Jakson, Classical Electrodynamics, 2nd edition, Wiley, (1975).
\item{}[Ku] Igor Kukavica, On the Dissipative Scale for the Navier-Stokes Equations, {\it IUMJ}, {\bf Vol. 38}, {\bf No. 3}, (1999).
\item{}[LL] Lieb, Loss, Analysis, {\it Grad Studies in Math}, {\bf 14}, AMS, (1997).
\item{}[RG] Ruzmaikina, Grujic, On depletion of the vortex-stretching term in the 3-D Navier-Stokes equations, {\it CMP}, {\bf Vol. 247}, {\bf No. 3}, (2004), 601-611.
\end{itemize}

\end{document}